\documentclass[11pt]{article}

 \usepackage{latexsym}
 \usepackage{amsbsy}
 \usepackage{amssymb}
 \usepackage{amsmath}
 \input{amssym.def}
 \input{amssym}
 \usepackage{graphicx}
 \usepackage{epsfig}
 \usepackage[latin1]{inputenc}%needed for bibliography
 \usepackage{a4wide}

 \newtheorem{thm}{Theorem}[section]
 \newtheorem{defin}[thm]{Definition}
 \newtheorem{lem}[thm]{Lemma}
 \newtheorem{prop}[thm]{Proposition}
 \newtheorem{cor}[thm]{Corollary}
 \newtheorem{rem}[thm]{Remark}
 \newtheorem{ex}[thm]{Example}

 \newcommand{\bthm}{\begin{thm}}
 \newcommand{\ethm}{\end{thm}}
 \newcommand{\bd}{\begin{defin}}
 \newcommand{\ed}{\end{defin}}
 \newcommand{\blem}{\begin{lem}}
 \newcommand{\elem}{\end{lem}}
 \newcommand{\bcor}{\begin{cor}}
 \newcommand{\ecor}{\end{cor}}
 \newcommand{\bprop}{\begin{prop}}
 \newcommand{\eprop}{\end{prop}}
 \newcommand{\brem}{\begin{rem} \rm}
 \newcommand{\erem}{\end{rem}}
 \newcommand{\bex}{\begin{ex} \rm}
 \newcommand{\eex}{\end{ex}}

 \newcommand{\pr}{\noindent{\bf Proof. }}
 \newcommand{\ep}{\nolinebreak{\hspace*{\fill}$\Box$ \vspace*{0.25cm}}}

 \newcommand{\beq}{\begin{equation}}
 \newcommand{\eeq}{\end{equation}}

 \newcommand{\bea}{\begin{eqnarray}}
 \newcommand{\eea}{\end{eqnarray}}

 \newcommand{\beas}{\begin{eqnarray*}}
 \newcommand{\eeas}{\end{eqnarray*}}

 \newcommand{\beqs}{\begin{equation*}}
 \newcommand{\eeqs}{\end{equation*}}

 \newcommand{\bi}{\begin{itemize}}
 \newcommand{\ei}{\end{itemize}}

 \newcommand{\ben}{\begin{enumerate}}
 \newcommand{\een}{\end{enumerate}}

 \newcommand{\ba}{\begin{array}}
 \newcommand{\ea}{\end{array}}

 \newcommand{\R}{\mathbb R}
 \newcommand{\N}{\mathbb N}
 \newcommand{\C}{\mathbb C}

 \newcommand{\cF}{\ensuremath{{\cal F}}}

 \newcommand{\cL}{\ensuremath{{\cal L}}}

 \newcommand{\cS}{\ensuremath{{\cal S}}}

 \newcommand{\eps}{\varepsilon}
 \newcommand{\vphi}{\varphi}

 \newcommand{\supp}{\mathop{\rm supp}\nolimits}

 \newcommand{\rep}{\mathop{\rm Re}\nolimits}
 \newcommand{\imp}{\mathop{\rm Im}\nolimits}
 \newcommand{\res}{\mathop{\rm Res}\nolimits}

 \newcommand{\tg}{\mathop{\rm tg}\nolimits}
 \newcommand{\ctg}{\mathop{\rm ctg}\nolimits}

 \newcommand{\tgh}{\mathop{\rm tgh}\nolimits}

 \newcommand{\notmid}{\mid\kern-0.5em\not\kern0.5em}

 \def\dj{d\kern-0.4em\char"16\kern-0.1em}
 \def\Dj{\mbox{\raise0.3ex\hbox{-}\kern-0.4em D}}

\begin{document}

 \title{Complex order fractional derivatives in viscoelasticity}

 \author{Teodor M. Atanackovi\'c
         \footnote{Faculty of Technical Sciences, Institute of Mechanics, University of Novi Sad,
         Trg D. Obradovi\'ca 6, 21000 Novi Sad, Serbia,
         Electronic mail: atanackovic@uns.ac.rs}\\
         Sanja Konjik
         \footnote{Faculty of Sciences, Department of Mathematics and Informatics, University of Novi Sad,
         Trg D. Obradovi\'ca 4, 21000 Novi Sad, Serbia,
         Electronic mail: sanja.konjik@dmi.uns.ac.rs}\\
         Stevan Pilipovi\'c
         \footnote{Faculty of Sciences, Department of Mathematics and Informatics, University of Novi Sad,
         Trg D. Obradovi\'ca 4, 21000 Novi Sad, Serbia,
         Electronic mail: pilipovic@dmi.uns.ac.rs}\\
         Du\v san Zorica
         \footnote{Institute of Mathematics, Serbian Academy of Sciences and Arts,
         Kneza Mihaila 36, 11000 Belgrade, Serbia,
         Electronic mail: dusan\_zorica@mi.sanu.ac.rs}
       }

 \date{}
 \maketitle

 \begin{abstract}
 We introduce complex order fractional derivatives in
 models that describe viscoelastic materials.
 This can not be carried out unrestrictedly, and therefore
 we derive, for the first time, real valued compatibility constraints,
 as well as physical constraints that lead to acceptable models.
 As a result, we introduce a new form of complex order fractional derivative.
 Also, we consider a fractional differential equation with complex derivatives,
 and study its solvability.
 Results obtained for stress relaxation and creep are illustrated by several
 numerical examples.

 \end{abstract}

%%%%%%%%%%%%%%%%%%%%%%%%%%%%%%%%%%%%%%%%%%%%%%%%%%
 \section{Introduction}
 \label{sec:intro}
%%%%%%%%%%%%%%%%%%%%%%%%%%%%%%%%%%%%%%%%%%%%%%%%%%

 Fractional calculus is a powerful tool for modeling various phenomena in
 mechanics, physics, biology, chemistry, medicine, economy, etc. Last few
 decades have brought a rapid expansion of the non-integer order
 differential and integral calculus, from which both the theory and its
 applications benefit significantly. However, most of the work done in this field
 so far has been based on the use of real order fractional derivatives and
 integrals. It is worth to mention that there are several authors who also
 applied complex order fractional derivatives to model various phenomena,
 see the work of Machado or Makris, \cite{Machado, Makris, MakrisConstantinou}.
 In all of these papers, restrictions on constitutive parameters that
 follow from the Second Law of Thermodynamics were not examined. In the
 analysis that follows this issue will be addressed.

 The main goal of this paper is to motivate and explain basic concepts
 of fractional calculus with complex order fractional derivatives.
 Throughout the paper we will investigate constitutive equations in the
 dimensionless form, for all independent ($t$ and $x$) and dependent
 ($\sigma$ and $\eps$) variables.
 Thus, consider a constitutive equation, given by \eqref{cce-gen}, connecting
 the stress $\sigma(t,x)$ at the point $x\in\R$ and time $t\in\R_+$, with
 the strain $\eps(t,x)$:
 \beq \label{cce-gen}
 \sum_{n=0}^{N}a_{n}\, {}_{0}D_{t}^{\alpha_{n}}\sigma(t,x)
 =\sum_{m=0}^{M}b_{m}\,{}_{0}D_{t}^{\beta_{m}}\eps(t,x),
 \eeq
 that contains fractional derivatives of complex order
 $\alpha_1,\ldots,\alpha_N,\beta_1,\ldots,\beta_M$.
 The precise definition of the operator  ${}_{0}D_{t}^{\eta}$
 of fractional differentiation with respect to $t$ is given below.
 In order to make a useful framework for the study of \eqref{cce-gen}
 we involve two types of conditions:
 1. real valued compatibility constraints, and
 2. thermodynamical constraints.
 Since this paper deals
 only with the well-posedness of constitutive equations of type \eqref{cce-gen}
 and their solvability for strain if stress is prescribed, we may, without loss
 of generality, assume that both $\sigma$ and $\eps$ are functions only of $t$.
 Also, equation \eqref{cce-gen} can be seen as a generalization of different
 models considered in the literature so far (see e.g.
 \cite{APSZ-book-vol1, CaputoMainardi-linmod, CaputoMainardi-dissmod, GonsovskiRossikhin, Hanyga-stf, Hanyga-tf, MainardiPagniniGorenflo, RS01-method}),
 since by taking all $\alpha_n$ and $\beta_m$ to be real numbers,
 the problem is reduced to the real case studied in the mentioned papers.

 Our results will show that constitutive equations of the form
 \eqref{cce-gen} may lead to creep and stress relaxation curves
 that are not monotonic.
 We note that the conditions of complete monotonicity required
 in e.g. \cite{AmendolaFabrizioGolden, FM, Mainardi}
 turn out not to be necessary but rather
 sufficient to meet the thermodynamic restrictions.
 It seems that conditions of monotonicity are stronger than the
 conditions following from the Second Law of Thermodynamics
 in the form of \cite{BagleyTorvik} analogously as in
 \cite{FL} where it is shown that the asymptotic stability,
 for certain class of constitutive equations, requires more
 extensive conditions on the coefficients compared with the
 restrictions that the classical formulation of the Second Law of
 Thermodynamics imposes.
 Non-monotonic creep curves were observed
 experimentally, however such a behavior is attributed to inertia
 either of the rod itself or of the rheometer.

 The paper is organized as follows. In Section \ref{sec:lfcoce} we
 investigate conditions leading to constitutive equations containing
 complex derivatives of stress and strain that can be used in viscoelastic
 models of the wave equation. More precisely, we derive restrictions
 on parameters in constitutive equation of the form \eqref{cce-gen} under which
 the Laplace and Fourier transforms, as well as their inverses, of real-valued functions will
 remain real-valued. Also, in order to impose the validity of the Second Law of
 Thermodynamics we follow the procedure presented in
 \cite{BagleyTorvik}, and obtain additional restrictions
 on model parameters in constitutive equations.
 As a result, we shall introduce a new form of the fractional
 derivative of complex order. Then, in Section \ref{sec:solvability},
 we treat the fractional Kelvin-Voigt complex order constitutive
 equation for viscoelastic body, find thermodynamical restrictions
 and prove sufficient conditions for its invertibility.
 Several numerical examples are presented in Section \ref{sec:numerics},
 as an illustration of creep and stress relaxation in viscoelastic materials.

 To the end of this section we recall basic definitions and
 results that will be used in our work.
 Fractional operators of complex order are introduced as follows
 (see \cite{Love-complex, SKM}):
 For $\eta\in\C$ with $0<\rep\eta<1$, definition of the left
 Riemann-Liouville fractional integral of an absolutely
 continuous function on $[0,T]$, $T>0$ ($y\in AC([0,T])$)
 coincides with the case of real $\eta$, i.e.,
 ${}_{0}I_{t}^{\eta}y(t) :=\frac{1}{\Gamma (\eta)}
 \int_{0}^{t} \frac{y(\tau)}{(t-\tau)^{1-\eta}}
 \,d\tau$, $t\in [0,T]$,
 where $\Gamma $ is the Euler gamma function.
 If $\eta=i\theta$, $\theta\in\R$, then the latter integral diverges,
 and hence one introduces the fractional integration of imaginary order as
 $$
 {}_{0}I_{t}^{i\theta}y(t) := \frac{d}{dt} {}_{0}I_{t}^{1+i\theta}y(t)
 = \frac{1}{\Gamma (1+i\theta)}
 \frac{d}{dt}\int_{0}^{t} (t-\tau)^{i\theta} y(\tau)\,d\tau,
 \quad t\in [0,T].
 $$
 However, in both cases the left Riemann-Liouville fractional derivative
 of order $\eta\in\C$ with $0\leq\rep\eta<1$ is given by
 $$
 {}_0D_t^{\eta} y(t) : = \frac{d}{dt} {}_{0}I_{t}^{1-\eta}y(t)
 = \frac{1}{\Gamma (1-\eta)}
 \frac{d}{dt}\int_{0}^{t} \frac{y(\tau)}{(t-\tau)^{\eta}}\,d\tau,
 \quad t\in [0,T].
 $$

 The basic tool for our study will be the Laplace and Fourier transforms.
 In order to have a good framework we will perform these transforms
 in $\cS'(\R)$, the space of tempered distributions. It is the dual space
 for the Schwartz space of rapidly decreasing functions $\cS(\R)$.
 In particular, we are interested in the space $\cS'_{+}(\R)$ whose
 elements are of the form $y=P(D)Y_0$, where $Y_0$ is a locally integrable
 polynomial bounded function on $\R$ that vanishes on $(-\infty,0)$,
 and $P(D)$ denotes a partial differential operator.

 The Fourier transform of $y\in L^1(\R)$ (or $y\in L^2(\R)$)
 is defined as
 \beq \label{ft}
 \cF y (\omega)=\hat{y}(\omega)=\int_{-\infty }^{\infty}
 e^{-i\omega x}y(x)\,dx, \qquad \omega \in \R.
 \eeq
 In the distributional setting, one has $\langle \cF y,\vphi \rangle
 =\langle y,\cF \vphi \rangle$, $y\in\cS'(\R)$, $\vphi \in \cS(\R)$,
 where $\cF\vphi$ is defined by \eqref{ft}.
 For $y\in L^{1}(\R)$ with $y(t)=0$, $t<0$, and $|y(t)| \leq Ae^{at}$, $a,A>0$,
 the Laplace transform is given by
 $$
 \cL y(s)= \tilde{y}(s) = \int_{0}^{\infty} e^{-st} y(t)\, dt,
 \qquad \rep s>a.
 $$
 If $y\in \cS'_{+}(\R)$ then $a=0$ (since $y$ is bounded by a polynomial).
 Then $\cL y$ is a holomorphic function in the half
 plane $\rep s>0$ (see e.g. \cite{Vladimirov-EqMP}).

 Let $Y(s)$, $\rep s>0$, be a holomorphic function bounded by a polynomial
 in that domain. Then, for a suitable polynomial $P$, $Y(s)/P(s)$ is integrable
 along the line $\Gamma=(a-i\infty, a+i\infty)$, and the inverse Laplace
 transform of $Y$ is a tempered distribution $y(t)=P(\frac{d}{dt})Y_0(t)$,
 where $Y_0(t) = \cL^{-1}[Y](t)=\frac{1}{2\pi i} \int_\Gamma \frac{Y(s)}{P(s)} e^{st}\, ds$.

 Let $y\in\cS'_+$. Recall:
 $$
 \cF\bigg[\frac{d^{n}}{dx^{n}}y\bigg] (\omega) = (i\omega)^{n}\cF y(\omega)
 \,\,\, (\omega\in\R),
 \quad
 \cL\bigg[\frac{d^{n}}{dt^{n}}y\bigg] (s) = s^{n}\cL y(s)
 \,\,\, (\rep s>0),
 \quad n\in\N,
 $$
 $$
 \cF[{}_{0}D_{x}^{\alpha}y](\omega) = (i\omega)^{\alpha}\cF y(\omega)
 \,\,\, (\omega\in\R),
 \quad
 \cL[{}_{0}D_{t}^{\alpha}y](s) = s^{\alpha}\cL y(s)
 \,\,\, (\rep s>0),
 \quad \alpha\in\C.
 $$

%%%%%%%%%%%%%%%%%%%%%%%%%%%%%%%%%%%%%%%%%%%%%%%%%%
 \section{Linear fractional constitutive equations with complex derivatives}
 \label{sec:lfcoce}
%%%%%%%%%%%%%%%%%%%%%%%%%%%%%%%%%%%%%%%%%%%%%%%%%%

 In what follows, we shall denote by $\alpha$ and $\beta$ the orders of
 fractional derivatives. $\alpha$ will be assumed to be a real number,
 while $\beta$ will be an element of $\C$ which is not real, i.e.,
 $\beta=A+iB$, with $B\not= 0$. Also, we shall assume that $0<\alpha, A<1$.

%%%%%%%%%%%%%%%%%%%%%%%%%%%%%%%%%%%%%%%%%%%%%%%%%%
 \subsection{Real valued compatibility constraints}
 \label{ssec:mc}
%%%%%%%%%%%%%%%%%%%%%%%%%%%%%%%%%%%%%%%%%%%%%%%%%%

 Similarly as in the real case, see \cite{AKOZ-thermodynamical},
 it is quite difficult to begin with a study of the most general case of \eqref{cce-gen}.
 Therefore, in order to try to recognize the essence of the problem
 and find possibilities for overcoming it,
 we shall first concentrate to simpler forms of constitutive equations that contain
 complex derivatives.
 Consider the following generalization of the Hooke law in the complex setting:
 \beq \label{cce-simplest}
 \sigma(t) = b\, {}_0D_t^{\beta} \eps(t),
 \eeq
 where $\beta\in\C$ and $b\in\R$. In order to find restrictions
 on parameters $b$ and $\beta$ in \eqref{cce-simplest} which
 yield a physically acceptable constitutive equation, we shall verify
 the next two conditions: For real strain $\eps$, the stress $\sigma$ has to be
 real valued function of $t$. We call this
 a real valued compatibility requirement.
 Thermodynamical restrictions will result from the Second
 Law of Thermodynamics, and will be studied in the next section.
 Note that in the case of constitutive equations with only real-valued
 fractional derivatives, the real valued compatibility requirement always holds true,
 while the thermodynamical restrictions had to be investigated
 (cf.\ \cite{AKOZ-thermodynamical}).

 \bthm \label{th:rvcc-1}
 Let $\eps\in AC([0,T])$ be real-valued, for all $T>0$, $0<A<1$ and $b\not=0$.
 Then function $\sigma$ defined by \eqref{cce-simplest} is real valued
 if and only if $\beta\in\R$.
 \ethm

 \pr
 It follows from \eqref{cce-simplest}, with $\beta=A+iB$,
 and $1/\Gamma(1-\beta)=h+ir$, that
 $$
 \imp\sigma(t) = b\,\frac{d}{dt} \int_0^t \eps(t-\tau) \tau^{-A}
 \Big(r\cos(B\ln\tau) - h\sin(B\ln\tau)\Big) \,d\tau,
 \quad t\geq 0.
 $$
 Denoting by $\frac{r}{h}:=\tg\phi$, for $h\not=0$, we obtain
 $$
 \imp\sigma(t) = \frac{bh}{\cos\phi} \,\frac{d}{dt} \int_0^t \eps(t-\tau) \tau^{-A}
 \sin(\phi-B\ln\tau)\,d\tau,
 \quad t\geq 0.
 $$
 In the case $h=0$, the imaginary part of $\sigma$ reduces to
 $br\frac{d}{dt} \int_0^t \eps(t-\tau) \tau^{-A} \cos(B\ln\tau)\,d\tau$.

 If $B=0$ then $r=0$ and $\phi=0$, hence $\imp\sigma=0$ and
 $\sigma$ is a real-valued function.

 Next, suppose that $B\not=0$. But then one can find a subinterval $(t_1,t_2)$
 of $(0,T)$ where $\sin(\phi-B\ln\tau)$, respectively $\cos(B\ln\tau)$,
 is positive (or negative), and choose $\eps\in AC([0,T])$ which is
 compactly supported in $(t_1,t_2)$ and strictly positive. This
 leads to a contradiction with the assumption $\imp\sigma=0$.
 \ep

 The previous theorem implies that equations of form \eqref{cce-simplest}
 with $\beta\in\C\backslash\R$ can not be a constitutive equation for
 a viscoelastic body.

 Next, consider the equation
 \beq \label{cce-eps}
 \sigma(t) = b_1\, {}_0D_t^{\beta_1} \eps(t)+ b_2\, {}_0D_t^{\beta_2} \eps(t),
 \qquad t\geq 0,
 \eeq
 where $b_1,b_2\in\R$, and $\beta_1,\beta_2\in\C\backslash\R$,
 i.e., $\beta_k=A_k+iB_k$ and $B_k\not= 0$ ($k=1,2$).
 Suppose again that $\eps\in AC([0,T])$ is a real-valued
 function, for every $T>0$.

 \brem
 Note that dimension $[{}_0D_t^{\beta_1} \eps]$ is $T^{-\beta_1}$,
 where $T$ is time unit. Therefore, \eqref{cce-eps} makes sense if
 $[b_1]=T^{\beta_1}$ and $[b_2]=T^{\beta_2}$.
 \erem

 \bthm \label{th:sufc}
 Function $\sigma$ given by \eqref{cce-eps} is real-valued for all
 real-valued positive $\eps\in AC([0,T])$ if and only if
 $b_1=b_2$ and $\beta_2=\bar{\beta_1}$.
 \ethm

 \pr
 We continue with the notation of Theorem \ref{th:rvcc-1}.
 Let $t\geq 0$. Then
 $$
 \sigma(t) =
 \frac{b_1}{\Gamma(1-\beta_1)} \frac{d}{dt} \int_0^t \eps(t-\tau) \tau^{-\beta_1}\, d\tau
 +
 \frac{b_2}{\Gamma(1-\beta_2)} \frac{d}{dt} \int_0^t \eps(t-\tau) \tau^{-\beta_2}\, d\tau
 $$
 Denote by $h_k+ir_k:=1/\Gamma(1-\beta_k)$, $k=1,2$. Then the imaginary part
 of the right hand side reads:
 \bea
 \imp\sigma(t) &=&
 \frac{d}{dt} \int_0^t \eps(t-\tau)
 \tau^{-A_1} \Big(b_1 r_1\cos(B_1\ln\tau) - b_1 h_1\sin(B_1\ln\tau)\Big)\,d\tau \nonumber \\
 && \,\, +
 \frac{d}{dt} \int_0^t \eps(t-\tau)
 \tau^{-A_2} \Big(b_2 r_2\cos(B_2\ln\tau) - b_2 h_2\sin(B_2\ln\tau)\Big)\,d\tau \label{sigma-im}
 \eea

 Using the identity $\Gamma(\bar{z})=\overline{\Gamma(z)}$,
 it is straight forward to check that  $b_1=b_2$ and $\beta_2=\bar{\beta_1}$
 imply that  $\imp\sigma=0$, and hence $\sigma$ is a
 real-valued function on $[0,T]$, for every $T>0$.

 Conversely, we want to find conditions on parameters
 which yield a real-valued function $\sigma$.
 Thus, we look at \eqref{sigma-im},
 with the change of variables $p=\ln\tau$, $\tau\in(0,t), t\leq T,$
 and solutions of the equation
 \bea
 && \frac{d}{dt} \int_{-\infty}^{\ln t} e^{p} \eps(t-e^p)
 \Big(b_1r_1 e^{-A_1 p} \cos(B_1p) - b_1h_1 e^{-A_1 p} \sin(B_1p) \nonumber \\
 && \qquad \qquad \qquad
 + b_2r_2 e^{-A_2 p} \cos(B_2 p) - b_2h_2 e^{-A_2 p} \sin(B_2 p) \Big)\,dp
 =0, \quad t\in[0,T]. \label{sig-im=0}
 \eea
 Set $\frac{r_k}{h_k}:=\tg \phi_k$, $k=1,2$, for $h_1,h_2\not= 0$.
 (For $h_k=0$ set $\phi_k:=\frac{\pi}{2}$, $k=1,2$.)
 Then \eqref{sig-im=0} gives
 $$
 \frac{d}{dt} \int_{-\infty}^{\ln t} e^{p} \eps(t-e^p) \bigg(
 \frac{b_1h_1}{\cos\phi_1} e^{-A_1 p} \sin(\phi_1-B_1p)
 + \frac{b_2h_2}{\cos\phi_2} e^{-A_2 p} \sin(\phi_2-B_2p)
 \bigg)\,dp = 0.
 $$

 Assume first that $|B_1|\not= |B_2|$, say $|B_1|>|B_2|$.
 Since the basic period of $\sin(\phi_1-B_1 p)$ ($T_0=2\pi/|B_1|$)
 is smaller than for $\sin(\phi_2-B_2 p)$, it follows that for every $k\in\N$,
 $k>k_0$, where $k_0$ depends on $\ln t$, the function $\sin(\phi_1-B_1p)$
 changes its sign at least three times in the interval $(k\pi,k\pi+2\pi/|B_2|)$.
 Thus, on that interval, there exist at least two intervals where
 $\sin(\phi_1-B_1p)$ and $\sin(\phi_2-B_2 p)$ have the same sign, and two
 intervals where they have opposite signs. We conclude that there exists
 an interval $[a,b]\subseteq (k\pi,k\pi+2\pi/|B_2|)\subset (-\infty, \ln t)$,
 so that
 \beq \label{iu}
 \frac{b_1h_1}{\cos\phi_1} e^{-A_1 p} \sin(\phi_1-B_1p)
 + \frac{b_2h_2}{\cos\phi_2} e^{-A_2 p} \sin(\phi_2-B_2p)
 >0, \quad p\in [a,b].
 \eeq
 Choose $ \delta>0$ and  $k\in\N$ so that
 $$
 \{t-e^p \,;\,  t\in(T/2-\delta,T/2+\delta), p\in[a,b] \}=(T/2-\delta-e^b,T/2+\delta-e^a)=I\subset(0,T).
 $$
 Now, we choose a non-negative function $\eps\in AC([0,T])$
 with the following properties: $\supp\eps\subseteq I$, so that
 the function $p\mapsto \eps(t-e^p)$, $p\in [a,b]$, is strictly positive on some
 $[a_1,b_1]\subseteq (a,b)$. This implies that for $t\in (T/2-\delta,T/2+\delta)$,
 $$
 \int_{a}^{b} e^p \eps(t-e^p) \bigg(
 \frac{b_1h_1}{\cos\phi_1} e^{-A_1 p} \sin(\phi_1-B_1p)
 + \frac{b_2h_2}{\cos\phi_2} e^{-A_2 p} \sin(\phi_2-B_2p)
 \bigg)\,dp
 $$
 is not a constant function. This is in contradiction with \eqref{sig-im=0}.

 Therefore, in order to have \eqref{sig-im=0}, one must have $|B_1|=|B_2|$.
 Moreover, arguing as above, one concludes that $|\phi_1-B_1p|=|\phi_2-B_2p|$
 must hold, for all $p\in(-\infty,\ln t)$, $t\geq 0$.
 Then we examine the equation
 $$
 b_1h_1 e^{-A_1 p}  - b_2h_2 e^{-A_2 p}= 0,
 \quad \mbox{ or } \quad
 b_1h_1 e^{-A_1 p}  + b_2h_2 e^{-A_2 p}= 0,
 \quad p\in[k\pi,k\pi+2\pi/|B_1|].
 $$
 %Because of analyticity, this equality holds for $p\in(-\ln t, \infty)$.
 Now in both cases, $b_1b_2>0$ or $b_1b_2<0,$ it is easy to conclude that
 $A_1=A_2$, $b_1=b_2$ and $B_1=-B_2$ have to be satisfied.
 This  proves the theorem.
 \ep

 \brem
 (i)\ Theorem \ref{th:sufc} states that a real valued compatibility constraint for
 constitutive equations of form \eqref{cce-eps}
 holds if they contain complex fractional derivatives of strain, whose
 orders have to be complex conjugated numbers.
 Therefore, we may assume in the sequel, without loss of generality,
 that $B>0$.

 (ii)\ According to the above analysis, one can take arbitrary linear
 combination of pairs of complex conjugated fractional derivatives of strain.
 Moreover, one can also allow the same type of fractional derivatives of stress.
 Thus, one can consider the most general stress-strain constitutive equation with
 fractional derivatives of complex order:
 $$
 \sigma(t) + \sum_{i=1}^{N} c_i \Big({}_0D_t^{\gamma_i} + {}_0D_t^{\bar{\gamma}_i}\Big) \sigma(t)
 =
 \eps(t) + \sum_{j=1}^{M} b_j \Big({}_0D_t^{\beta_j} + {}_0D_t^{\bar{\beta}_j}\Big) \eps(t),
 $$
 where $c_i,b_j\in\R$ and $\gamma_i,\beta_j\in\C$, $i=1,\ldots,N$, $j=1,\ldots,M$.

 (iii)\ As a consequence one has that stress-strain relations can also
 contain arbitrary real order fractional derivatives, without any additional
 restrictions.
 This fact has already been known from previous work.

 (iv)\ The same conclusions can also be obtained using a different approach.
 One can apply the result from  \cite[p.\ 293, Satz 2]{Doetsch-Handbuch-1},
 which tells that a function $F$ is real-valued (almost everywhere), if its
 Laplace transform is real-valued, for all real $s$ in the half-plane of
 convergence right from some real $x_0$, in order to show that an
 admissible fractional constitutive equation \eqref{cce-gen} may be of
 complex order only if it contains pairs of complex conjugated fractional
 derivatives of stress and strain.
 \erem

%%%%%%%%%%%%%%%%%%%%%%%%%%%%%%%%%%%%%%%%%%%%%%%%%%
 \subsection{Thermodynamical restrictions}
 \label{ssec:tr-complex}
%%%%%%%%%%%%%%%%%%%%%%%%%%%%%%%%%%%%%%%%%%%%%%%%%%

 In the analysis that follows we shall consider the isothermal
 processes only. For such processes the Second Law of
 Thermodynamics, i.e., the entropy of the system increases,
 is equivalent to the dissipativity condition
 \beq \label{disseq}
 \int_0^T \sigma(t) \dot{\eps}(t)\, dt \geq 0.
 \eeq
 Equation \eqref{disseq} is one-dimensional version
 of the Second Principle of Thermodynamics for simple
 materials under isothermal conditions
 (see \cite[p.\ 83]{AmendolaFabrizioGolden}).
 Also, in writing \eqref{disseq} we observe the fact that
 for $t\in(-\infty,0)$ the system is in virginal state.
 Note that $T=\infty$ in \cite[p.\ 113]{Day}.
 Since we consider arbitrary $T$, condition \eqref{disseq}
 is stronger than the condition in \cite{Day}.

 In \cite{AmendolaFabrizioGolden} stress is assumed as
 \beq \label{AFG-1}
 \sigma(t)=G_0 \eps(t) + \int_0^\infty G'(u)\eps(t-u)\,du,
 \eeq
 where $G_0$ is instanteneous elastic modulus,
 and $G$ is relaxation function,
 while stain is taken in the form
 $$
 \eps(s) = \eps_0\cos(\omega s) + \eps_0 \sin(\omega s),
 \qquad s\leq t,\,\,\, \omega>0,
 $$
 where $\eps_0$ is the amplitude of strain.
 These assumptions along with \eqref{disseq} and $T=\frac{2n\pi}{\omega}$
 lead to
 \beq \label{AFG-2}
 \int_0^\infty G'(u) \sin(\omega u) \, du < 0.
 \eeq
 Note that \eqref{AFG-2} is equivalent to
 $$
 \cF_s[G'(t)](\omega) = \imp\Big(\cF[G'(t)](\omega)\Big) > 0,
 $$
 where $\cF_s[f(t)](\omega)=\int_0^\infty f(t)\sin(\omega t)\,
 dt$.

 Since constitutive equations containing fractional derivatives
 that we shall treat in this work do not obey \eqref{AFG-1},
 we follow the approach of \cite{BagleyTorvik}.
 Namely, in \cite{BagleyTorvik} it was assumed that periodic
 strain $\eps$ results in periodic stress $\sigma$ after
 the end of transient regime. Positivity of dissipation
 work \eqref{disseq} for a cycle and each time instant
 within the cycle leads to
 $$
 \rep \hat{E}(\omega) \geq 0 \quad \mbox{ and } \quad
 \imp \hat{E}(\omega) \geq 0, \qquad \forall\,\omega>0,
 $$
 where $\hat{E}=\frac{\hat{\sigma}}{\hat{\eps}}$ is the
 complex modulus obtained from the constitutive
 equation by applying the Fourier transform
 (cf.\ \cite[Eq. (20), (21)]{BagleyTorvik}).
 Note that $\rep \hat{E}$ and $\imp \hat{E}$ are referred as
 storage and loss modulus respectively .

 We start with the constitutive equation
 \beq \label{cce}
 \sigma(t) = 2b\, {}_0\bar{D}_t^\beta \eps(t),
 \qquad
 {}_0\bar{D}_t^\beta := \frac12 \Big({}_0D_t^{\beta} + {}_0D_t^{\bar{\beta}}\Big),
 \qquad t\geq 0,
 \eeq
 where we assume that $b>0$ and $\beta=A+iB$, $0<A<1$, $B>0$.
 Note that \eqref{cce} generalizes the Hooke law in the complex fractional framework.
 In the case $\beta\in\R$, this new complex fractional operator
 ${}_0\bar{D}_t^\beta$ coincides with the usual left Riemann-Liouville
 fractional derivative.

 We apply the Fourier transform to \eqref{cce}:
 $\hat{\sigma}(\omega) =
 b\, ((i\omega)^{\beta} + (i\omega)^{\bar{\beta}}) \hat{\eps}(\omega)$,
 $\omega\in \R$.
 Then, define the complex modulus $\hat{E}$ such that
 $\hat{\sigma}(\omega) = \hat{E}(\omega) \cdot \hat{\eps}(\omega)$,
 $\omega\in\R$, i.e.,
 $$
 \hat{E}(\omega) : = b\, \big((i\omega)^{\beta} + (i\omega)^{\bar{\beta}}\big)
 = b\, \omega^A \big(e^{-\frac{B\pi}{2}} e^{i(\frac{A\pi}{2} + B\ln\omega)}
 +e^{\frac{B\pi}{2}} e^{i(\frac{A\pi}{2} - B\ln\omega)}\big),
 \qquad \omega \in \mathbb R.
 $$
 Thermodynamical restrictions involve, for $\omega\in\R_+$,
 \bea
 \rep \hat{E}(\omega) &=& b\, \omega^A
 \Big(
 e^{-\frac{B\pi}{2}} \cos\Big(\frac{A\pi}{2} + B\ln\omega\Big)
 + e^{\frac{B\pi}{2}} \cos\Big(\frac{A\pi}{2} - B\ln\omega\Big)
 \Big) \geq 0, \label{tdr-re} \\
 \imp \hat{E}(\omega) &=& b\, \omega^A
 \Big(
 e^{-\frac{B\pi}{2}} \sin\Big(\frac{A\pi}{2} + B\ln\omega\Big)
 + e^{\frac{B\pi}{2}} \sin\Big(\frac{A\pi}{2} - B\ln\omega\Big)
 \Big) \geq 0. \label{tdr-im}
 \eea
 But this is in contradiction with $B>0$, since for $\omega>0$,
 \eqref{tdr-re} and \eqref{tdr-im} imply $B=0$.

 In order not to confront the real valued compatibility requirement
 and the Second Law of Thermodynamics for \eqref{cce},
 one may require that \eqref{tdr-re} and \eqref{tdr-im} hold for $\omega$ in some
 bounded interval instead of in all of $\R$.
 Alternatively, as we shall do in the sequel, one can
 modify \eqref{cce} by adding additional terms,  in order to preserve
 the Second Law of Thermodynamics.

 Thus, we proceed by proposing the following constitutive equation
 \beq \label{rcce}
 \sigma(t) = a\, {}_0D_t^\alpha \eps(t) + 2b\, {}_0\bar{D}_t^\beta \eps(t),
 \qquad t\geq 0,
 \eeq
 where we assume that $a,b>0$, $\alpha\in\R$, $0<\alpha<1$, and $\beta=A+iB$, $0<A<1$, $B>0$.
 Again we follow the procedure described above for
 deriving thermodynamical restrictions:
 $\hat{\sigma}(\omega) = [a(i\omega)^\alpha +
 b ((i\omega)^{\beta} + (i\omega)^{\bar{\beta}})]
 \hat{\eps}(\omega)$, $\omega\in \R$.
 Consider the complex module ($\omega\in\R$)
 \beq \label{cmpl-mdl}
 \hat{E}(\omega) = a\,\omega^\alpha e^{i\frac{\alpha\pi}{2}} +
 b\, \omega^A \Big(e^{-\frac{B\pi}{2}} e^{i(\frac{A\pi}{2} + B\ln\omega)}
 +e^{\frac{B\pi}{2}} e^{i(\frac{A\pi}{2} - B\ln\omega)}\Big);
 \eeq
 \bea
 \rep \hat{E}(\omega) &=& a\, \omega^\alpha \cos\frac{\alpha\pi}{2}
 + b\, \omega^A
 \Big(
 e^{-\frac{B\pi}{2}} \cos \Big(\frac{A\pi}{2} + B\ln\omega\Big)
 + e^{\frac{B\pi}{2}} \cos \Big(\frac{A\pi}{2} - B\ln\omega\Big)
 \Big), \label{re-n} \\
 \imp \hat{E}(\omega) &=& a\,\omega^\alpha \sin\frac{\alpha\pi}{2}
 + b\, \omega^A
 \Big(
 e^{-\frac{B\pi}{2}} \sin \Big(\frac{A\pi}{2} + B\ln\omega\Big)
 + e^{\frac{B\pi}{2}} \sin \Big(\frac{A\pi}{2} - B\ln\omega\Big)
 \Big). \label{im-n}
 \eea
 We will investigate conditions $\rep\hat{E}\geq 0$ and $\imp\hat{E}\geq 0$ on
 $\R_+$. The assumption $\alpha>A$ leads to a contradiction, since for
 $\omega\searrow 0$ the sign of the second term in \eqref{re-n} determines
 the sign of $\rep\hat{E}$, and it can be negative. Thus, we must have $\alpha\leq A$.
 If $\alpha<A$. then for $\omega\to\infty$,
 the second term in \eqref{re-n} could be negative. This together yields that
 the only possibility is $A=\alpha$.
 (The same conclusion is obtained if one considers $\imp\hat{E}\geq 0$, $\omega>0$.)

 Therefore, with $A=\alpha$, \eqref{re-n} and \eqref{im-n} become:
 \bea
 \rep \hat{E}(\omega) &=& a\, \omega^\alpha \cos\frac{\alpha\pi}{2}
 + 2b\, \omega^\alpha f(\omega), \qquad \omega>0, \label{re-Aal} \\
 \imp \hat{E}(\omega) &=& a\, \omega^\alpha \sin\frac{\alpha\pi}{2}
 + 2b\, \omega^\alpha g(\omega), \qquad \omega>0, \label{im-Aal}
 \eea
 with
 \bea
 f(\omega) & := & \cos\frac{\alpha\pi}{2}\cos(\ln\omega^B)\cosh\frac{B\pi}{2}
 + \sin\frac{\alpha\pi}{2}\sin(\ln\omega^B)\sinh\frac{B\pi}{2},
 \qquad \omega>0, \label{f-g-f} \\
 g(\omega) & := & \sin\frac{\alpha\pi}{2}\cos(\ln\omega^B)\cosh\frac{B\pi}{2}
 - \cos\frac{\alpha\pi}{2}\sin(\ln\omega^B)\sinh\frac{B\pi}{2},
 \qquad \omega>0. \label{f-g-g}
 \eea
 We further have that $\rep\hat{E}(\omega)\geq a\omega^\alpha\cos\frac{\alpha\pi}{2}
 +2b\omega^\alpha\min_{\omega\in\R_+} f(\omega)$, $\omega>0$, and the similar estimate for
 $\imp\hat{E}$, hence we now look for the minimums of functions $f$ and $g$ on $\R_+$.
 Using the substitution $x=\ln\omega^B$ we find that $f'(x)=0$ and $g'(x)=0$
 at $x_f$ and $x_g$ so that
 $$
 \tg x_f = \tg\frac{\alpha\pi}{2} \tgh\frac{B\pi}{2}
 \qquad \mbox{ and } \qquad
 \tg x_g = - \ctg\frac{\alpha\pi}{2} \tgh\frac{B\pi}{2}.
 $$
 Solutions $x_{f1}, x_{f2}, x_{g1}$ and $x_{g2}$ satisfy:
 $x_{f1}\in (0,\frac{\pi}{2})$, $x_{f2}\in (\pi ,\frac{3\pi}{2})$,
 and $x_{g1}\in (\frac{\pi}{2},\pi)$, $x_{g2}\in (\frac{3\pi}{2},2\pi)$,
 since $\tg \frac{\alpha \pi}{2}, \ctg \frac{\alpha \pi}{2}, \tgh \frac{B\pi}{2}>0$,
 and
 \begin{eqnarray*}
 f(x_{f}) &=& \pm \cos\frac{\alpha \pi}{2} \cosh \frac{B\pi}{2}
 \sqrt{1+\Big(\tg \frac{\alpha \pi}{2} \tgh \frac{B\pi}{2}\Big)^{2}}, \\
 g(x_{g}) &=& \pm \sin\frac{\alpha \pi}{2} \cosh \frac{B\pi}{2}
 \sqrt{1+\Big(\ctg \frac{\alpha \pi}{2} \tgh \frac{B\pi}{2}\Big)^{2}}.
 \end{eqnarray*}
 Therefore, we have
 $\min_{x\in \R} f(x) = f(x_{f2})$ and $\min_{x\in \R}g(x) = g(x_{g1})$,
 so that \eqref{re-Aal} and \eqref{im-Aal} can be estimated as
 \beas
 \rep\hat{E}(\omega) & \geq & \omega^\alpha \cos\frac{\alpha \pi}{2} \bigg(
 a-2b \cosh \frac{B\pi}{2}\sqrt{1+\Big(\tg \frac{\alpha \pi}{2} \tgh \frac{B\pi}{2}\Big)^{2}}
 \bigg), \qquad \omega>0, \\
 \imp\hat{E}(\omega) & \geq & \omega^\alpha\sin\frac{\alpha \pi}{2} \bigg(
 a-2b\cosh \frac{B\pi}{2}\sqrt{1+\Big(\ctg \frac{\alpha \pi}{2} \tgh \frac{B\pi}{2}\Big)^{2}}
 \bigg), \qquad \omega>0.
 \eeas
 We obtain the thermodynamical restrictions for \eqref{rcce} by requiring
 $\rep\hat{E}(\omega) \geq 0$ and $\imp\hat{E}(\omega) \geq 0$,
 for $\omega \in \R_{+}$:
 \beq \label{rstr-0-12}
 a \geq 2b \cosh \frac{B\pi}{2}\sqrt{1+\Big(\ctg \frac{\alpha \pi}{2} \tgh \frac{B\pi}{2}\Big)^{2}},
 \qquad \mbox{ if } \,\,\, \alpha \in \Big(0,\frac{1}{2}\Big],
 \eeq
 and
 \beq \label{rstr-12-1}
 a \geq 2b \cosh \frac{B\pi}{2}\sqrt{1+\Big(\tg \frac{\alpha \pi}{2}\tgh \frac{B\pi}{2}\Big)^{2}},
 \qquad \mbox{ if } \,\,\, \alpha \in \Big[\frac{1}{2},1\Big).
 \eeq
 Notice that both restrictions further imply $a\geq 2b$.

 \brem
 (i)\ Fix $a$ and $\alpha$. Inequalities \eqref{rstr-0-12} and \eqref{rstr-12-1}
 impliy that as $B$ increases, the constant $b$ has to decrease,
 i.e., the contribution of complex fractional derivative of strain in the
 constitutive equation \eqref{rcce} is smaller if its imaginary part is larger.

 (ii)\ Also, inequalities \eqref{rstr-0-12} and \eqref{rstr-12-1} lead to the same
 restrictions on parameters $a, b, \alpha$ and $B$, since for $\alpha\in(0,\frac12]$
 one has $1-\alpha\in[\frac12,1)$, and the values of \eqref{rstr-0-12} and \eqref{rstr-12-1}
 coincide.

 (iii)\ Under the same conditions, constitutive equation \eqref{rcce} can be extended to
 $\sigma(t) = \eps(t) + a\, {}_0D_t^\alpha \eps(t) + 2b\, {}_0\bar{D}_t^\beta \eps(t)$,
 $t\geq 0$, which will be investigated in the next section.
 \erem

%%%%%%%%%%%%%%%%%%%%%%%%%%%%%%%%%%%%%%%%%%%%%%%%%%
 \section{Complex order fractional Kelvin-Voigt model}
 \label{sec:solvability}
%%%%%%%%%%%%%%%%%%%%%%%%%%%%%%%%%%%%%%%%%%%%%%%%%%

 Consider the constitutive equation involving the complex order fractional derivative
 \beq \label{nkj}
 \sigma(t) = \eps(t) + a\, {}_0D_t^\alpha \eps(t)  + 2b\, {}_0\bar{D}_t^\beta \eps(t),
 \qquad t\geq 0.
 \eeq
 We assume $a,b, E>0$, $0<\alpha<1$, $B>0$, $\beta=\alpha+iB$,
 and $\sigma$ and $\eps$ are real-valued functions.
 Note that \eqref{nkj} is a generalization of the model proposed
 in \cite{RS01-dugo}. Namely, \eqref{nkj} agrees with the latter
 when $\beta$ is real and positive.
 The inverse relation, i.e., $\eps$ as a function of $\sigma$ is given in Theorem \ref{th:lst}.

 The Laplace transform of \eqref{nkj} is
 $\tilde{\sigma}(s) = E \big(1+a\,s^\alpha+b\,(s^\beta+s^{\bar{\beta}})\big) \tilde{\eps}(s)$,
 $\rep s>0$, and hence
 \beq \label{eps-tilda}
 \tilde{\eps}(s) = \frac{1}{1+a\,s^\alpha+b\,(s^\beta+s^{\bar{\beta}})} \tilde{\sigma}(s),
 \qquad \rep s>0.
 \eeq
 In order to determine $\eps$ from \eqref{eps-tilda} we need to analyze zeros of
 \begin{equation} \label{psi}
 \psi(s) = 1+ a\,s^\alpha+b\, \big(s^\beta+s^{\bar{\beta}} \big), \quad s\in\C.
 \end{equation}
 Note that if we put $s=i\omega$, $\omega \in \R_+$, in \eqref{psi},
 it becomes the complex modulus:
 \begin{equation} \label{comp-mod}
 \psi(i\omega) = 1+ \hat{E}(\omega) =
 1+a\,(i\omega)^{\alpha}+b\,\big((i\omega)^{\beta}+(i\omega)^{\bar{\beta}}\big),
 \qquad \omega>0,
 \end{equation}
 where $\hat{E}$ is given in \eqref{cmpl-mdl}.

 Let $s=re^{i\vphi}$, $r>0$, $\vphi\in[0,2\pi]$. Then (with $\beta =\alpha +iB$)
 $$
 \psi (s) = 1+ar^\alpha e^{i\alpha \vphi}+br^\alpha\big(
 e^{-B\vphi } e^{i(\ln r^{B}+\alpha \vphi)}+ e^{B\vphi } e^{-i(\ln r^{B}-\alpha \vphi)}
 \big),
 $$
 and
 \beq \label{Re-1}
 \rep\psi (s) =1+ar^\alpha \cos(\alpha \vphi) + 2br^\alpha \big(
 \cos(\ln r^{B}) \cos(\alpha\vphi) \cosh(B\vphi) + \sin(\ln r^{B})
 \sin(\alpha \vphi) \sinh(B\vphi)
 \big),
 \eeq
 \beq \label{Im-1}
 \imp\psi (s) = ar^\alpha \sin(\alpha \vphi) + 2br^\alpha \big(
 \cos(\ln r^{B}) \sin(\alpha \vphi) \cosh(B\vphi)
 - \sin(\ln r^{B}) \cos(\alpha \vphi) \sinh(B\vphi)
 \big).
 \eeq

%%%%%%%%%%%%%%%%%%%%%%%%%%%%%%%%%%%%%%%%%%%%%%%%%%
 \subsection{Thermodynamical restrictions}
 \label{ssec:tr}
%%%%%%%%%%%%%%%%%%%%%%%%%%%%%%%%%%%%%%%%%%%%%%%%%%

 In the case of \eqref{comp-mod}, using \eqref{Re-1} and \eqref{Im-1} we obtain:
 \beas
 \rep\psi(i\omega) &=& 1+\rep \hat{E}(\omega)
 \geq 1+a\omega^\alpha \cos\frac{\alpha \pi}{2}
 + 2b\omega^\alpha \min_{x\in \R} f(x),
 \qquad x=\ln\omega^B, \, \omega>0, \\
 \imp\psi(i\omega) &=& \imp \hat{E}(\omega)
 \geq a\omega^\alpha \sin\frac{\alpha \pi}{2}
 +2b\omega^\alpha \min_{x\in \R} g(x),
 \qquad x=\ln\omega^B, \, \omega>0,
 \eeas
 where $f$ and $g$ are as in \eqref{f-g-f} and \eqref{f-g-g}.
 This leads to the same thermodynamical restrictions \eqref{rstr-0-12} and
 \eqref{rstr-12-1}, as in Section \ref{ssec:tr-complex}.

 Therefore, from now on, we shall assume \eqref{rstr-0-12} and \eqref{rstr-12-1}
 to hold true. Now we shall examine the zeros of $\psi$.

%%%%%%%%%%%%%%%%%%%%%%%%%%%%%%%%%%%%%%%%%%%%%%%%%%
 \subsection{Zeros of $\psi$ and solutions of \eqref{nkj}}
 \label{ssec:zeros}
%%%%%%%%%%%%%%%%%%%%%%%%%%%%%%%%%%%%%%%%%%%%%%%%%%

 \bthm \label{th:zeros}
 Let $\psi$ be the function given by \eqref{psi}. Then
 \bi
 \item[(i)] $\psi$ has no zeros in the right complex half-plane $\rep s\geq 0$.
% if the coefficients $a, b, \alpha$ and $B$ satisfy thermodynamical restrictions
% \eqref{rstr-0-12} and \eqref{rstr-12-1}.
 \item[(ii)] $\psi$ has no zeros in $\C$ if the coefficients $a, b, \alpha$ and $B$ satisfy
% \eqref{rstr-0-12}, \eqref{rstr-12-1}, and
 \beq \label{bs-1}
 \begin{split}
 a & \geq 2b \cosh (B\pi) \sqrt{1+ \big(\tg(\alpha \pi) \tgh(B\pi) \big)^{2}},
 \quad \mbox{ for } \,\,\, \alpha\in\Big[\frac14,\frac34\Big]\backslash\Big\{\frac12\Big\}, \\
 a & \geq 2b \cosh (B\pi) \sqrt{1+ \big(\ctg(\alpha \pi) \tgh(B\pi) \big)^{2}}.
 \quad \mbox{ for } \,\,\, \alpha\in\Big(0,\frac14\Big) \cup \Big\{\frac12\Big\} \cup \Big(\frac34, 1\Big).
 \end{split}
 \eeq
 \ei
 \ethm

 \pr
 First, we notice that if $s_{0}$ is a solution to $\psi (s) = 0$,
 then $\bar{s}_{0}$ (the complex conjugate of $s_{0}$) is also a solution,
 since $\psi(\bar{s}) = 1+a\,\bar{s}^\alpha +b\,(\bar{s}^\beta + \bar{s}^{\bar{\beta}})
 = \overline{\psi(s)}$.
 Thus, we restrict our attention to the upper complex half-plane $\imp s\geq0$,
 i.e., $\vphi \in [0,\pi]$.

 Using (\ref{Re-1}) and (\ref{Im-1}) we have, with $s=re^{i\vphi}$, $r>0$, $\vphi\in[0,\pi]$,
 and $x=\ln r^B$,
 \begin{eqnarray}
 \rep \psi (s) &\geq &1+ar^{\alpha}\cos (\alpha \vphi) +2br^{\alpha} \min_{x\in\R} f(x),  \label{min-ef-1} \\
 \imp \psi (s) &\geq & ar^{\alpha}\sin (\alpha \vphi) +2br^{\alpha} \min_{x\in \R}g(x),  \label{min-ge-1}
 \end{eqnarray}
 where
 \begin{eqnarray}
 f(x) &=&\cos (x) \cos (\alpha \vphi) \cosh (B \vphi)
 + \sin (x) \sin (\alpha\vphi) \sinh (B \vphi),
 \qquad x\in\R,  \label{ef-1} \\
 g(x) &=& \cos (x) \sin (\alpha \vphi) \cosh (B \vphi)
 - \sin (x) \cos (\alpha \vphi) \sinh (B \vphi),
 \qquad x\in\R.  \label{ge-1}
 \end{eqnarray}
 The critical points $x_{f}$ and $x_{g}$ of $f$ and $g$, respectively, satisfy
 \beq \label{kriticni}
 \tg x_{f} = \tg (\alpha \vphi) \tgh (B \vphi) \geq 0
 \quad \mbox{ and } \quad
 \tg x_{g}=-\ctg (\alpha \vphi) \tgh (B \vphi) \leq 0,
 \eeq

 The proof of (i) and (ii) will be given by the argument principle.

\begin{figure}[htbp]
 \centering
 \includegraphics[width=7.7cm]{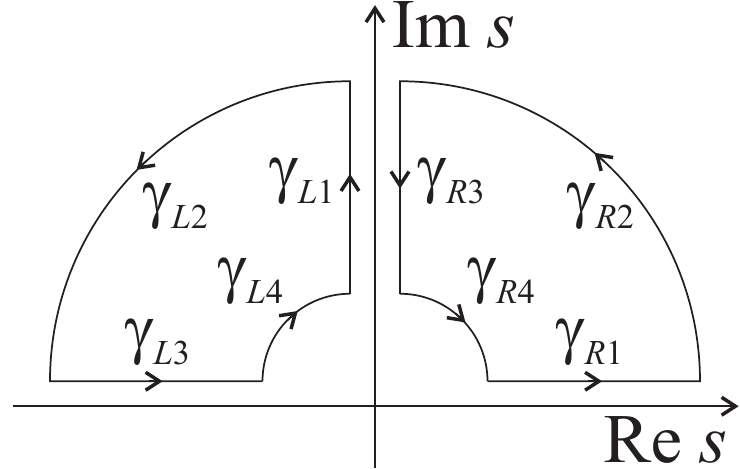}
 \caption{Contour $\Gamma$.}
 \label{fig-cont}
 \end{figure}

 (i)\ Consider $\psi$ in the case $\rep s, \imp s >0$. Choose a contour
 $\Gamma= \gamma_{R1}\cup\gamma_{R2}\cup\gamma_{R3}\cup\gamma_{R4}$,
 as it is shown in Fig. \ref{fig-cont}.

 $\gamma_{R1}$ is parametrized by $s=x$, $x\in(\eps,R)$ with $\eps\to 0$ and $R\to\infty$,
 so that \eqref{Re-1} and \eqref{Im-1} yield
 \begin{eqnarray*}
 \rep\psi (x) &=&1+x^{\alpha} (a+2b\cos (\ln x^{\beta})) \geq 1+x^{\alpha} (a-2b) \geq 0, \\
 \imp\psi (x) &=& 0,
 \end{eqnarray*}
 since both \eqref{rstr-0-12} and \eqref{rstr-12-1} imply $a\geq 2b$.
 Moreover, we have $\lim_{x\to 0} \psi(x)=1$
 and $\lim_{x\to \infty} \psi(x)=\infty$.

 Along $\gamma_{R2}$ one has $s=Re^{i\vphi}$, $\vphi \in [ 0,\frac{\pi}{2}] $,
 with $R\rightarrow \infty$. By \eqref{kriticni} we have $\tg x_{f} \geq 0$, so
 that
 $$
 \min_{x\in \R}f(x) = -\cos (\alpha \vphi) \cosh (B \vphi) \sqrt{1+( \tg (\alpha \vphi) \tgh(B \vphi) )^{2}},
 $$
 and therefore \eqref{min-ef-1} becomes
 $$
 \rep \psi (s) \geq 1+R^{\alpha}\cos (\alpha \vphi) (a-2b\cosh (B\vphi)
 \sqrt{1+( \tg(\alpha \vphi) \tgh (B \vphi) )^{2}}) \geq 0.
 $$
 The previous inequality holds true, since for $\vphi \in [ 0,\frac{\pi}{2}] $ we have that
 $$
 p(\vphi) = \cosh (B \vphi) \sqrt{1+(\tg (\alpha \vphi) \tgh (B \vphi))^{2}}
 \leq \cosh \frac{B \pi}{2}\sqrt{1+\Big(\tg \frac{\alpha\pi}{2}\tgh \frac{B \pi}{2}\Big)^{2}},
 $$
 because of the fact that the function $p$ monotonically increases
 on $[0,\frac{\pi}{2}]$. Moreover, by \eqref{Re-1} and \eqref{Im-1}, we have
 \begin{eqnarray*}
 \rep \psi (s) &\rightarrow & \infty  \quad \mbox{ and } \quad \imp\psi (s) = 0,
 \quad \mbox{ for } \quad
 \vphi=0, \,\,\, R\rightarrow \infty, \\
 \rep \psi (s) &\rightarrow & \infty,
 \quad \mbox{ for } \quad
 \vphi =\frac{\pi}{2}, \,\,\, R\rightarrow \infty.
 \end{eqnarray*}

 The next segment is $\gamma_{R3}$, which is parametrized by
 $s=i\omega$, $\omega\in[R,\eps]$, with $\eps\rightarrow 0$ and $R\rightarrow \infty$.
 Then \eqref{Re-1} and \eqref{Im-1} yield
 $$
 \rep \psi (i\omega ) = 1 + \rep\hat{E} (\omega) \geq 0
 \quad \mbox{ and } \quad
 \imp \psi (i\omega) = \imp\hat{E}  (\omega) \geq 0,
 \,\,\, \omega \in (\eps,R),
 $$
 due to the thermodynamical requirements. Moreover, by \eqref{Re-1} and \eqref{Im-1},
 we have
 \begin{eqnarray*}
 \rep\psi (\omega) &\rightarrow &1
 \quad \mbox{ and } \quad
 \imp\psi (\omega) \rightarrow 0,
 \quad \mbox{ as } \quad
 \omega \rightarrow 0, \\
 \rep\psi (\omega) &\rightarrow &\infty
 \quad \mbox{ and } \quad
 \imp\psi (\omega) \rightarrow \infty,
 \quad \mbox{ as } \quad
 \omega \rightarrow \infty.
 \end{eqnarray*}

 The last part of the contour $\Gamma$ is the arc $\gamma_{R4}$, with
 $s=\eps e^{i\vphi}$, $\vphi \in [0,\frac{\pi}{2}]$, with $\eps \rightarrow 0$.
 Using the same arguments as for the contour $\gamma_{R2}$, we have
 \begin{equation} \label{re-na-desnom-luku}
 \rep\psi (s) \geq 1+\eps^\alpha \cos (\alpha \vphi)
 \big(a-2b \cosh (B\vphi) \sqrt{1+(\tg (\alpha \vphi) \tgh (B\vphi))^{2}}\big) \geq 1.
 \end{equation}
 Also, by \eqref{Im-1} and \eqref{re-na-desnom-luku}, we have
 $$
 \rep\psi (s) \rightarrow 1
 \quad \mbox{ and } \quad
 \imp\psi (s) \rightarrow 0,
 \quad \mbox{ as } \quad
 \eps \rightarrow 0.
 $$

 We conclude that $\Delta \arg \psi (s)=0$ so that, by the argument principle,
 there are no zeroes of $\psi$ in the right complex half-plane $\rep s\geq 0$.

 (ii)\ In order to discuss the zeros of $\psi$ in the left complex half-plane,
 we use the contour $\Gamma_{L}=\gamma_{L1}\cup \gamma_{L2}\cup
 \gamma_{L3}\cup \gamma_{L4}$, shown in Fig. \ref{fig-cont}.
 The contour $\gamma_{L1}$ has the same parametrization as the
 contour $\gamma_{R3}$, so the same conclusions as for $\gamma_{R3}$ hold
 true.

 The parametrization of the contour $\gamma_{L2}$ is $s=Re^{i\vphi}$,
 $\vphi \in [\frac{\pi}{2},\pi]$, with $R\rightarrow \infty$.
 Let us distinguish two cases.

 \begin{enumerate}
 \item[a1.] $\alpha \vphi \in (0,\frac{\pi}{2})$ \\
 Then $\sin (\alpha\vphi) >0$, $\cos (\alpha\vphi) > 0$,
 so the critical points of $f$ and $g$, see \eqref{ef-1} and \eqref{ge-1},
 given by \eqref{kriticni}, satisfy $\tg x_{f}>0$ and $\tg x_{g}< 0$.
 For the minimums of $f$ and $g$ in \eqref{min-ef-1} and
 \eqref{min-ge-1} we have
 \begin{eqnarray*}
 \min_{x\in \R} f(x) &=& -\cos (\alpha\vphi) \cosh (B\vphi)
 \sqrt{1+ (\tg (\alpha\vphi) \tgh (B\vphi))^{2}}, \\
 \min_{x\in \R} g(x) &=& -\sin (\alpha\vphi) \cosh (B\vphi)
 \sqrt{1+ (\ctg (\alpha\vphi) \tgh (B\vphi))^{2}},
 \end{eqnarray*}
 respectively, so that \eqref{min-ef-1} and \eqref{min-ge-1} become
 \begin{eqnarray}
 \rep \psi (s) & \geq &1+R^\alpha \cos (\alpha \vphi)
 \Big(a-2b\cosh (B\vphi) \sqrt{1+\big(\tg(\alpha\vphi) \tgh (B\vphi) \big)^{2}} \Big), \label{a1-repsi} \\
 \imp \psi (s) & \geq & R^\alpha \sin (\alpha \vphi)
 \Big(a-2b\cosh (B\vphi) \sqrt{1+\big( \ctg(\alpha\vphi) \tgh (B\vphi) \big)^{2}}\Big). \label{a1-impsi}
 \end{eqnarray}
 Function $H_{f}(\vphi) = \cosh (B\vphi)
 \sqrt{1+(\tg (\alpha\vphi) \tgh (B\vphi))^{2}}$ is monotonically increasing
 for $\vphi \in [\frac{\pi}{2},\pi]$, since $\alpha \vphi \in (0,\frac{\pi}{2})$,
 thus
 $$
 \rep \psi (s) \geq 1+R^\alpha \cos (\alpha \vphi) \Big( a-2b \cosh (B\pi)
 \sqrt{1+\big( \tg (\alpha \pi) \tgh (B\pi))^{2}}\Big) \geq 0,
 $$
 if \eqref{bs-1} is satisfied. Note that $\rep \psi (s) \rightarrow \infty$
 and $\imp \psi (s) \rightarrow \infty$, for $\vphi = \pi$, $R\rightarrow \infty$.

 \item[b1.] $\alpha \vphi \in [\frac{\pi}{2},\pi)$ \\
 Then $\sin (\alpha\vphi) >0$, $\cos (\alpha\vphi)\leq 0$, so the critical
 points of $g$, see \eqref{ge-1}, given by \eqref{kriticni}, satisfy
 $\tg x_{g}\geq 0$. For the minimum of $g$ in \eqref{ge-1}
 we have
 $$
 \min_{x\in \R} g(x) = -\sin (\alpha\vphi) \cosh (B\vphi)
 \sqrt{1+\big( \ctg (\alpha\vphi) \tgh (B\vphi) \big)^{2}},
 $$
 so that \eqref{min-ge-1} becomes
 $$
 \imp \psi (s) \geq R^\alpha \sin (\alpha \vphi)
 \Big( a-2b \cosh (B\vphi) \sqrt{1+\big(\ctg(\alpha\vphi) \tgh (B\vphi) \big)^{2}}\Big).
 $$
 Function $H_{g}(\vphi) =\cosh (B\vphi)
 \sqrt{1+(\ctg (\alpha\vphi) \tgh (B\vphi))^{2}}$ is monotonically
 increasing for $\vphi \in [\frac{\pi}{2},\pi]$, since
 $\alpha \vphi \in [\frac{\pi}{2},\pi)$, and \eqref{bs-1} implies
 $$
 \imp \psi (s) \geq R^\alpha \sin (\alpha \vphi) \Big( a-2b \cosh (B\pi)
 \sqrt{1+\big( \ctg (\alpha \pi) \tgh (B\pi) \big)^{2}}\Big) \geq 0.
 $$
 Note that $\imp \psi (s) \rightarrow \infty$, for $\vphi =\pi$, $R\rightarrow \infty$.
 \end{enumerate}

 Now we discuss possible situations for $\alpha\in(0,1)$
 and $\vphi\in[\frac{\pi}{2},\pi]$.

 If $\alpha\in(0,\frac12)$ then a1. holds so that $\rep\psi(s)\geq 0$.

 If $\alpha\in[\frac12,1)$ then we distinguish two cases.
 For $\vphi\in[\frac{\pi}{2}, \frac{\pi}{2\alpha})$ case a1. holds, so
 $\rep\psi(s)\geq 0$.
 For $\vphi\in[\frac{\pi}{2\alpha}, \pi)$ case b1. holds, and $\imp\psi(s)\geq 0$.

% Let $\alpha \in (0,\frac{1}{2})$, and $\vphi \in [\frac{\pi}{2},\pi]$.
% Then a1. holds, so that $\rep\psi (s) \geq 0$ for all $s\in \gamma_{L2}$.
% On the other hand, for $\alpha \in [\frac{1}{2},1)$,
% there exists $\vphi_{0}\in (\frac{\pi}{2},\pi)$
% so that for $\vphi \in [\frac{\pi}{2},\vphi_{0})$ case a1. holds,
% i.e., $\rep \psi (s) \geq 0$, while for $\vphi \in (\vphi_{0},\pi]$
% case b1. holds, i.e., $\imp \psi (s) \geq 0$.

 Parametrization of the contour
 $\gamma_{L3}$ is $s=xe^{i\pi}$, $x\in (\eps ,R)$, with $\eps \rightarrow 0$
 and $R\rightarrow \infty$. Again, we have two cases.

 \begin{enumerate}
 \item[a2.] $\alpha \in (0,\frac{1}{2})$ \\
 Then, for $x\in (\eps ,R)$,
 using the same argumentation as in case a1, we have, by
 \eqref{a1-repsi} and \eqref{a1-impsi} (with $R=x$),
 $\rep\psi(s)\geq 1$ and $\imp\psi(s)\geq 0$,
 due to \eqref{bs-1}.
 Thus, looking at  \eqref{a1-repsi} and \eqref{a1-impsi} (with $R=x$),
 we conclude
 \begin{eqnarray*}
 \rep \psi (s) & \rightarrow &\infty
 \quad \mbox{ and } \quad
 \imp \psi (s) \rightarrow \infty,
 \quad \mbox{ for } \quad
 x \rightarrow \infty, \\
 \rep \psi (s) &\rightarrow &1
 \quad \mbox{ and } \quad
 \imp \psi (s) \rightarrow 0,
 \quad \mbox{ for } \quad
 x\rightarrow 0.
 \end{eqnarray*}

 \item[b2.] $\alpha \in [\frac{1}{2},1)$ \\
 Then, for $x\in (\eps ,R)$,
 using the same argumentation as in case b1, we have (by \eqref{a1-impsi})
 that $\imp\psi(s)\geq 0$, and
 $$
 \imp \psi (s) \rightarrow \infty,
 \quad \mbox{ for } \quad
 x\rightarrow \infty
 \quad \mbox{ and } \quad
 \imp \psi (s) \rightarrow 0,
 \quad \mbox{ for } \quad
 x \rightarrow 0.
 $$
 \end{enumerate}

 The parametrization of the contour $\gamma_{L4}$ is $s=\eps e^{i\vphi}$,
 $\vphi \in [ \frac{\pi}{2},\pi]$, with $\eps \rightarrow 0$.
 From \eqref{Re-1} and \eqref{Im-1}, for sufficiently small $\eps$, we have
 $$
 \rep \psi (s) \rightarrow 1
 \quad \mbox{ and } \quad
 \imp \psi (s) \rightarrow 0,
 \quad \mbox{ for } \quad
 \eps \rightarrow 0,\, \vphi \in \Big[ \frac{\pi}{2},\pi \Big].
 $$

 Summing up all results from cases a1, a2, b1 and b2,
 we obtain the following:
 \bi
 \item For $\alpha \in (0,\frac{1}{2})$,
 $\rep \psi (s) \geq 0$, for $s\in \Gamma_{L}$, which implies that
 $\Delta \arg \psi (s) =0$. Therefore, using the argument principle,
 we conclude that in this case $\psi$ has no zeros in the left complex half-plane.
 \item If $\alpha \in [\frac{1}{2},1)$, then for $s\in \gamma_{L1}$
 and $s\in \{z\in \gamma_{L2} \,|\, \arg z\leq\frac{\pi}{2\alpha}\}$,
 we have $\rep \psi (s) \geq 0$, while for
 $s\in \{z\in \gamma_{L2} \,|\, \arg z>\frac{\pi}{2\alpha}\}$ and $s\in \gamma_{L3}$,
 we have $\imp \psi (s) \geq 0$. For $s\in \gamma_{L4}$ we again have
 $\rep \psi (s) \geq 0$.
 Hence, we conclude that $\Delta \arg \psi (s) =0$, and therefore,
 using the argument principle, neither in case $\alpha \in [\frac{1}{2},1)$
 function $\psi$ has zeros in the left complex half-plane.
 \ei
 This completes the proof.
 \ep

 Rewrite \eqref{eps-tilda} as
 \beq \label{eps-k-tilda}
 \tilde{\eps}(s) = \tilde{K}(s) \tilde{\sigma}(s),
 \qquad
 \tilde{K}(s) = \frac{1}{1+as^{\alpha}+b(s^{\beta}+s^{\bar{\beta}})},
 \qquad \rep s>0.
 \eeq

 \bthm \label{th:lst}
 Let $\tilde{\eps}$ be given by \eqref{eps-k-tilda}. Then
 \begin{equation} \label{eps}
 \eps (t) = K(t) \ast \sigma (t), \quad t\geq 0.
 \end{equation}

 Moreover, if \eqref{bs-1} holds, then
 \bea
 K(t) = K_{I}(t) & = & \frac{1}{2\pi i} \int_{0}^{\infty} \Bigg(
 \frac{e^{-qt}}{1+q^{\alpha} e^{i\alpha \pi} \Big[
 a + b(e^{i\ln q^B} e^{-B \pi}+ e^{- i\ln q^B} e^{B \pi})
 \Big]} \nonumber \\
 & & \qquad -
 \frac{e^{-qt}}{1+q^{\alpha} e^{-i\alpha \pi} \Big[
 a + b(e^{i\ln q^B} e^{B \pi}+ e^{- i\ln q^B} e^{-B \pi})
 \Big]}
 \Bigg)\,dq.  \label{k}
 \eea

 If condition \eqref{bs-1} is violated, then $\psi$ has at most finite number
 of zeros in the left complex half-plane, and
 $$
 K=K_{I}
 \qquad \mbox{ or } \qquad
 K=K_{I}+K_{R},
 $$
 where
 \beq \label{K-R}
 K_{R}(t)  = \sum_{\psi(s_i)=0 \atop i=1,2,\ldots, n} \Big(
 \res(\tilde{K}(s)  e^{st},s_i) +\res(\tilde{K}(s)  e^{st},\bar{s}_i)
 \Big),
 \eeq
 with $\tilde{K}$ given by \eqref{eps-k-tilda}.
 \ethm

\begin{figure}[htbp]
 \centering
 \includegraphics[width=7.7cm]{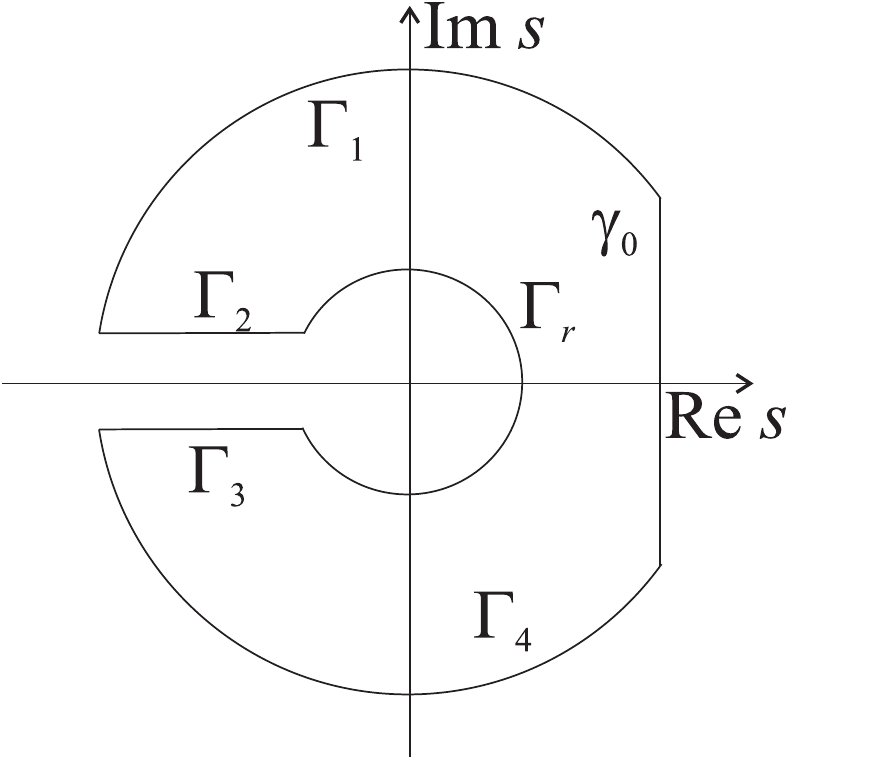}
 \caption{Contour $\tilde{\Gamma}$.}
 \label{hk}
 \end{figure}

 \pr
 The first part is clear. We invert now $\tilde{K}$, given in \eqref{eps-k-tilda},
 by the use of the Cauchy residues theorem
 \begin{equation} \label{crt}
 \oint_{\tilde{\Gamma}}\tilde{K}(s)  e^{st}\,ds=2\pi i\sum_{\psi(\tilde{s})=0}\res(\tilde{K}(s)  e^{st}, \tilde{s})
 \end{equation}
 and the contour $\tilde{\Gamma} =\Gamma _{1}\cup \Gamma _{2}\cup \Gamma _{r}\cup
 \Gamma _{3}\cup \Gamma _{4}\cup \gamma _{0}$ shown in Fig. \ref{hk}.

 If condition \eqref{bs-1} is satisfied,
 then, by Theorem \ref{th:zeros}, the residues equal zero.
 One can show that the integrals over the contours
 $\Gamma _{1}$, $\Gamma _{r}$ and $\Gamma _{4}$ tend to zero when
 $R\rightarrow \infty $ and $r\rightarrow 0$. The remaining integrals give:
 \begin{eqnarray*}
 \lim_{\substack{ R\rightarrow \infty, \\ r\rightarrow 0}} \int_{\Gamma _{2}} \tilde{K}(s)  e^{st}\,ds
 &=& - \int_{0}^{\infty} \frac{e^{-qt}}{1+q^{\alpha} e^{i\alpha \pi} \Big[
 a + b(e^{i\ln q^B} e^{-B \pi}+ e^{- i\ln q^B} e^{B \pi})
 \Big]} \,dq, \\
 \lim_{\substack{ R\rightarrow \infty, \\ r\rightarrow 0}}\int_{\Gamma _{3}} \tilde{K}(s)  e^{st}\,ds
 &=& \int_{0}^{\infty} \frac{e^{-qt}}{1+q^{\alpha} e^{-i\alpha \pi} \Big[
 a + b(e^{i\ln q^B} e^{B \pi}+ e^{- i\ln q^B} e^{-B \pi})
 \Big]} \,dq, \\
 \lim_{\substack{ R\rightarrow \infty, \\ r\rightarrow 0}}\int_{\gamma _{0}} \tilde{K}(s)  e^{st}\,ds
 &=&2\pi i K_{I}(t),
 \end{eqnarray*}
 which, by the Cauchy residues theorem \eqref{crt}, leads to \eqref{k}.

 If condition \eqref{bs-1} is violated,
 then, by Theorem \ref{th:zeros}, the denominator of $\tilde{K}$
 either has no zeros in the complex plane, and so $K=K_{I}$,
 or it has zeros in the left complex half-plane,
 which comes in pairs with complex conjugates.
 We show now that $\psi$ (which is the dominator of $\tilde{K}$)
 can have at most finite number of zeros for $\imp s\leq 0$.
 Rewrite $\psi(s)=0$ as $a+b\, (s^{iB} + s^{-iB}) = s^{-\alpha}$.
 If $s=re^{i\vphi}$, $\vphi\in[\frac{\pi}{2}, \pi]$, then we have
 $$
 a+b\, (r^{iB}e^{-\vphi B} + r^{-iB} e^{\vphi B}) = r^{-\alpha} e^{-\alpha \vphi i}.
 $$
 When $r\to\infty$ the right hand side tends to zero, while
 the left hand side tends to $a$. As $r\to 0$ we see that the left hand side
 is bounded, while the right hand side is not bounded. Thus, the zeros
 in the left half-plane of function $\psi$, if exist, have to be bounded both from
 above and below. In that case we have $K=K_{I}+K_{R}$, where $K_R$ is
 given by \eqref{K-R}.
 \ep
 
 \brem
 In the cases of standard and fractional linear solid model creep compliance
 and relaxation modulus are exponential functions and derivative of a one-parameter
 Mittag-Leffler function. Also, these models have different relaxation and retardation
 times as well as the limiting values of the material functions. For more complicated
 models, as our is, functions representing creep compliance and relaxation modulus
 are possibly not so well known special functions. We note that other generalizations
 of standard linear viscoelastic solids, e.g. the one presented in
 \cite{RS01-dugo},
 also do not tackle this question.
 \erem

%%%%%%%%%%%%%%%%%%%%%%%%%%%%%%%%%%%%%%%%%%%%%%%%%%
 \section{Numerical verifications}
 \label{sec:numerics}
%%%%%%%%%%%%%%%%%%%%%%%%%%%%%%%%%%%%%%%%%%%%%%%%%%

 Here we present several examples of the proposed constitutive equation.
 We shall treat stress relaxation, creep and periodic loading cases.

%%%%%%%%%%%%%%%%%%%%%%%%%%%%%%%%%%%%%%%%%%%%%%%%%%
 \subsection{Stress relaxation experiment}
 \label{ssec:stress rel}
%%%%%%%%%%%%%%%%%%%%%%%%%%%%%%%%%%%%%%%%%%%%%%%%%%

 We take \eqref{nkj} with $\eps(t) = H(t)$, $H$ is the
 Heaviside function, and regularize it as $H_{\eps}(t) =1-\exp (-t/k)$,
 $k\rightarrow 0$. In order to determine $\sigma$ we calculate
 \beq \label{n1}
 \sigma (t) = H_{\eps}(t) +a\,{}_{0} D_{t}^{\alpha} H_{\eps}(t)
 +b\,{}_{0} \bar{D}_{t}^{\beta} H_{\eps}(t),
 \quad t\geq 0,
 \eeq
 with $\sigma(0) =0$, using the expansion formula (see \cite{AS04-1, AS08-1}),
 \beq \label{n2}
 {}_{0}D_{t}^{\gamma} y(t) \approx \frac{y(t)}{t^{\gamma}} A(N,\gamma)
 - \sum_{p=1}^{N} C_{p-1}(\gamma) \frac{V_{p-1}(y)(t)}{t^{p+\gamma}},
 \eeq
 where
 $$
 A(N,\gamma) =
 \frac{\Gamma (N+1+\gamma)}{\gamma \Gamma(1-\gamma) \Gamma(\gamma) N!},
 \qquad
 C_{p-1}(\gamma) = \frac{\Gamma (p+\gamma)}{\Gamma(1-\gamma) \Gamma(\gamma)(p-1)!},
 $$
 and
 \begin{equation} \label{n3}
 V_{p-1}^{(1)}(y)(t) = t^{p-1}y(t),
 \quad V_{p-1}(y)(0) =0,
 \quad p=1,2,3,\ldots
 \end{equation}
 Inserting \eqref{n2} into \eqref{n1} we obtain
 \begin{eqnarray}
 \sigma (t) &\approx& \bigg\{
 1+\bigg[
 \frac{a\,A(N,\alpha)}{t^{\alpha}}
 +b\bigg(\frac{A(N,\alpha +iB)}{t^{\alpha +iB}}+
 \frac{A(N,\alpha -iB)}{t^{\alpha -iB}}\bigg)
 \bigg]
 \bigg\}
 H_{\eps}(t)  \nonumber \\
 && - \sum_{p=1}^{N} \bigg\{
 \frac{C_{p-1}(\alpha)}{t^{p+\alpha}}+\bigg[
 \frac{C_{p-1}(\alpha +iB)}{t^{p+\alpha +iB}}
 +\frac{C_{p-1}(\alpha - iB)}{t^{p+\alpha -iB}}
 \bigg]
 \bigg\}
 V_{p-1}(H_\eps)(t),  \label{n4}
 \end{eqnarray}
 where
 \begin{equation} \label{n5}
 V_{p-1}^{(1)}(H_\eps)(t)=t^{p-1}H_{\eps}(t) ,
 \quad p=1,2,3,\ldots
 \end{equation}

 We will compare \eqref{n4} with the stress $\sigma $ obtained by \eqref{n1}
 and by the definition of fractional derivative \eqref{cce}:
 \begin{eqnarray}
 \sigma (t) &=& H_{\eps}(t)+a\frac{d}{dt}\frac{1}{\Gamma (1-\alpha)}
 \int_{0}^{t} \frac{H_{\eps}(\tau) \,d\tau}{(t-\tau)^{\alpha}}  \nonumber \\
 && + b\,\frac{d}{dt}
 \bigg[\frac{1}{\Gamma (1-\alpha -iB)}
 \int_{0}^{t} \frac{H_{\eps}(\tau)\,d\tau}{(t-\tau)^{\alpha +iB}}
 +\frac{1}{\Gamma (1-\alpha +iB)}
 \int_{0}^{t}\frac{H_{\eps}(\tau)\,d\tau}{(t-\tau)^{\alpha -iB}}\bigg]. \label{s-N6}
 \end{eqnarray}

\begin{figure}[htbp]
 \centering
 \includegraphics[width=7.7cm]{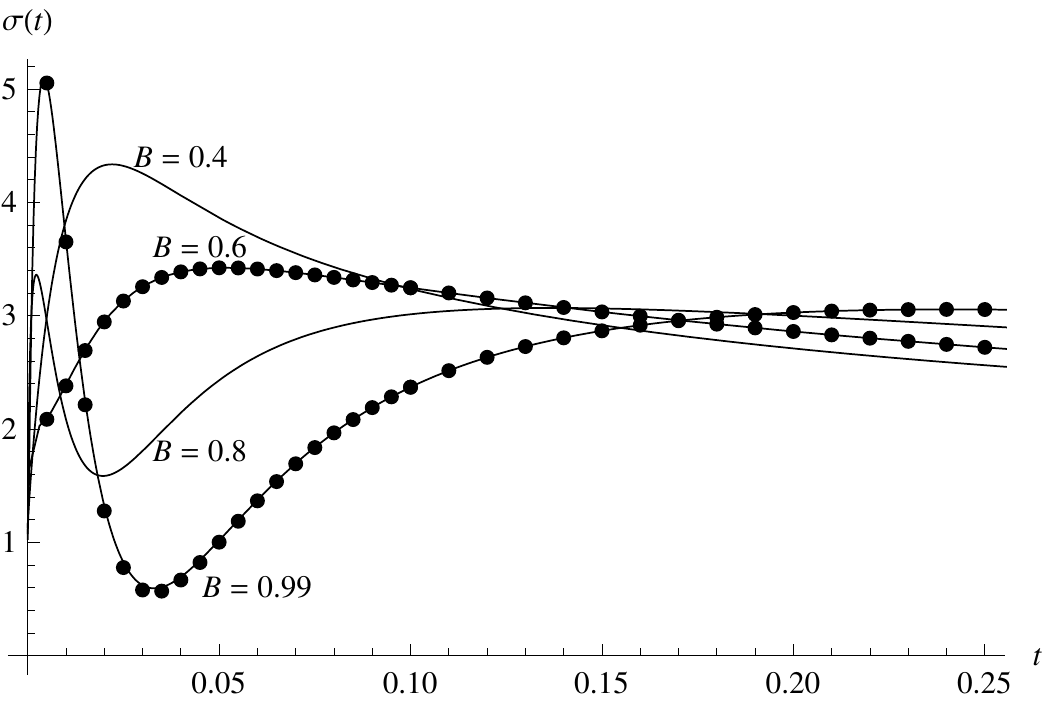}
 \caption{Stress relaxation curves for $\protect\alpha =0.4$ and
 $B\in \{0.4,0.6,0.8,0.99\}$, $a=0.8$, $b=0.1$, $t\in [0,0.25]$.}
 \label{fig-stressrel-1}
 \end{figure}

\begin{figure}[htbp]
 \centering
 \includegraphics[width=7.7cm]{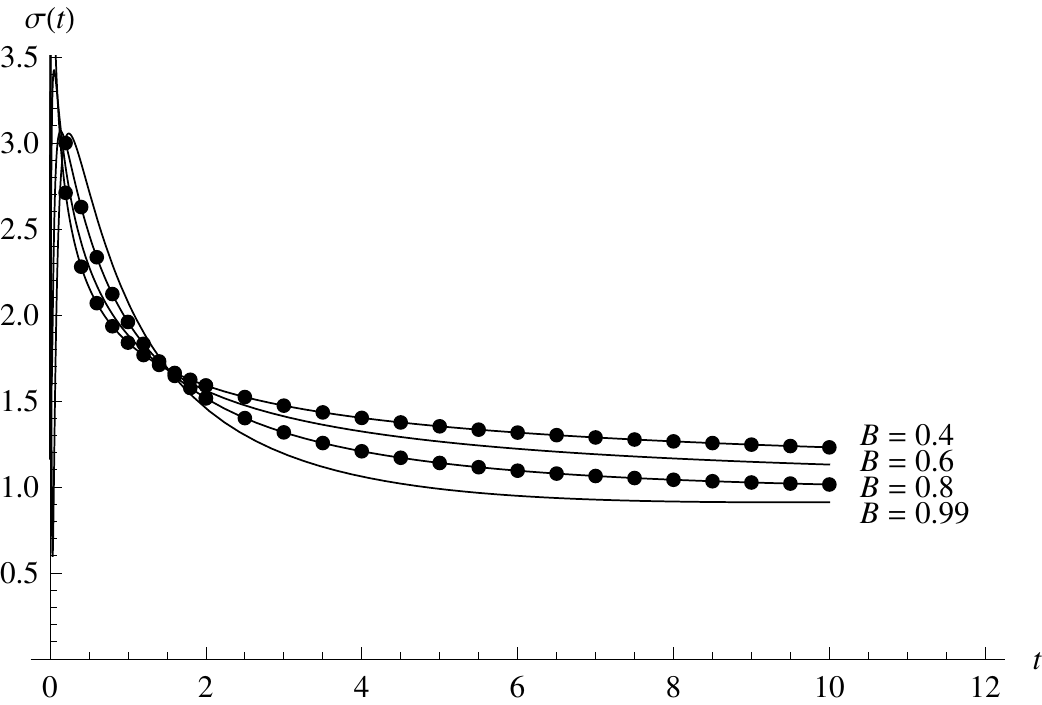}
 \caption{Stress relaxation curves for $\protect\alpha =0.4$ and
 $B\in \{0.4,0.6,0.8,0.99\}$, $a=0.8$, $b=0.1$, $t\in [ 0,10]$.}
 \label{fig-stressrel-2}
 \end{figure}

 In Fig. \ref{fig-stressrel-1} we show results obtained by determining
 $\sigma$ from \eqref{s-N6} for small times and different values of $B$. In the
 same figure we show, by dots, the values of $\sigma$, in several points,
 obtained by using \eqref{n4}, \eqref{n5} for $k=0.01$, $N=100$. As could be
 seen from Fig. \ref{fig-stressrel-1}, the results obtained from
 \eqref{s-N6} and \eqref{n4}, \eqref{n5} agree well.
 The stress relaxation curves are shown in Fig. \ref{fig-stressrel-2},
 for the same set of parameters and for
 larger times. As could be seen, regardless of the values of $B$, we have
 $\lim_{t\to \infty} \sigma (t) =1$. Note that in all cases of $B$, the restriction
 which follows from the dissipation inequality is satisfied.

%%%%%%%%%%%%%%%%%%%%%%%%%%%%%%%%%%%%%%%%%%%%%%%%%%
 \subsection{Creep experiment}
 \label{ssec:creep}
%%%%%%%%%%%%%%%%%%%%%%%%%%%%%%%%%%%%%%%%%%%%%%%%%%

\begin{figure}[htbp]
 \centering
 \includegraphics[width=7.7cm]{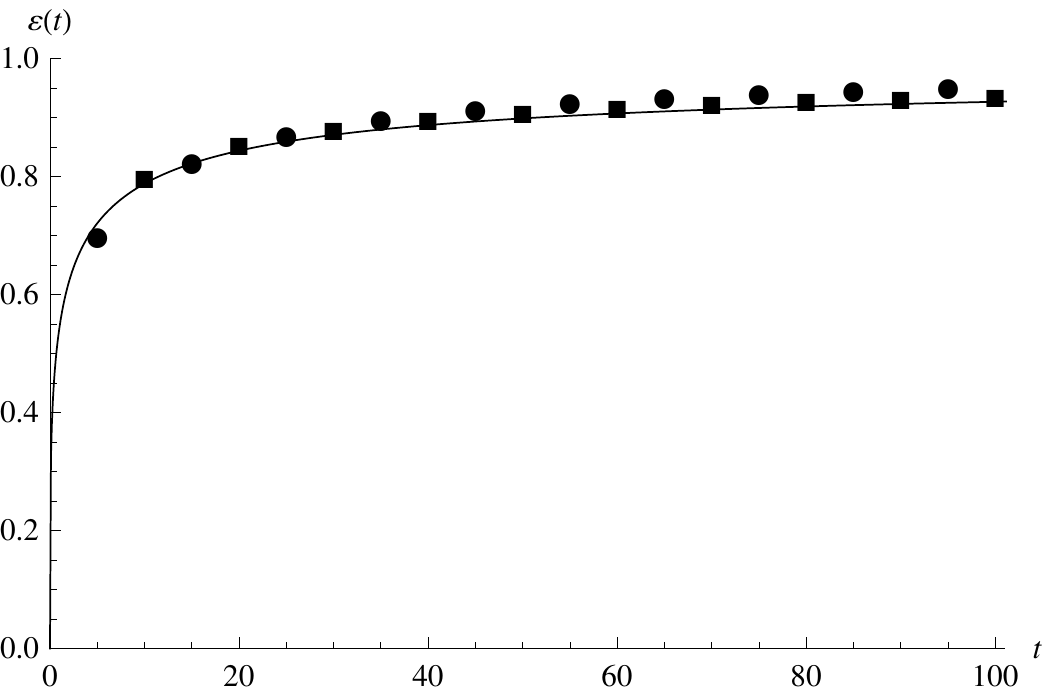}
 \caption{Creep curve for $\protect\alpha =0.4$ and $B=0.4$, $a=0.8$, $b=0.1$,
 $t\in [0,100]$.}
 \label{fig-creep-3}
 \end{figure}

\begin{figure}[htbp]
 \centering
 \includegraphics[width=7.7cm]{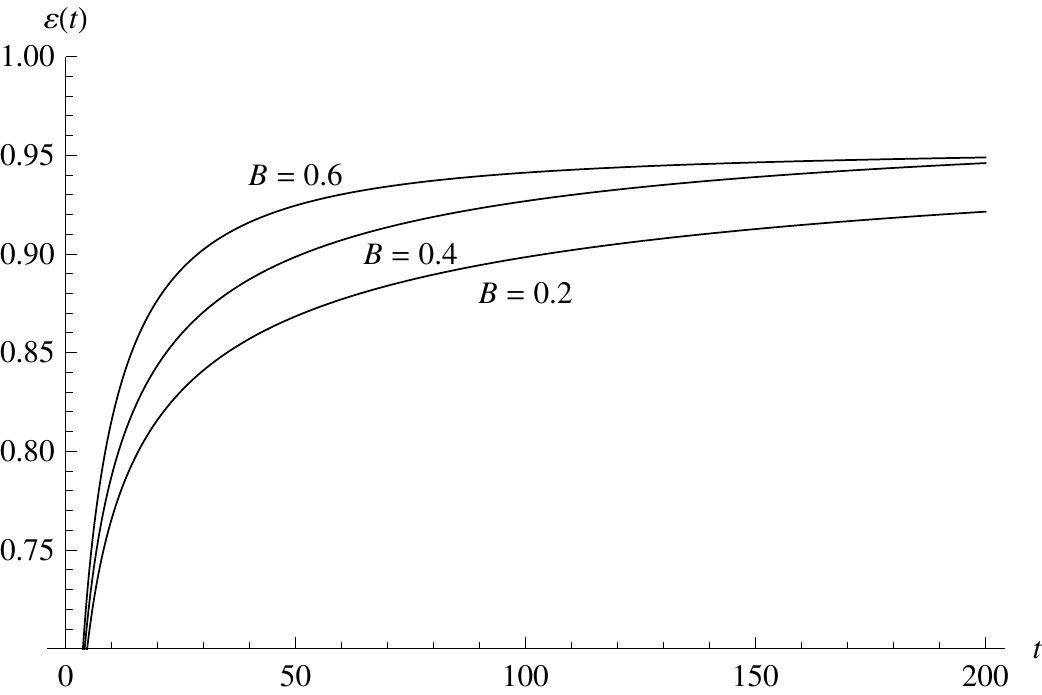}
 \caption{Creep curves for $\protect\alpha =0.4$ and $B\in \{0.2,0.4,0.6\}$,
 $a=0.8$, $b=0.1$, $t\in [0,200]$.}
 \label{fig-creep-1}
 \end{figure}

\begin{figure}[htbp]
 \centering
 \includegraphics[width=7.7cm]{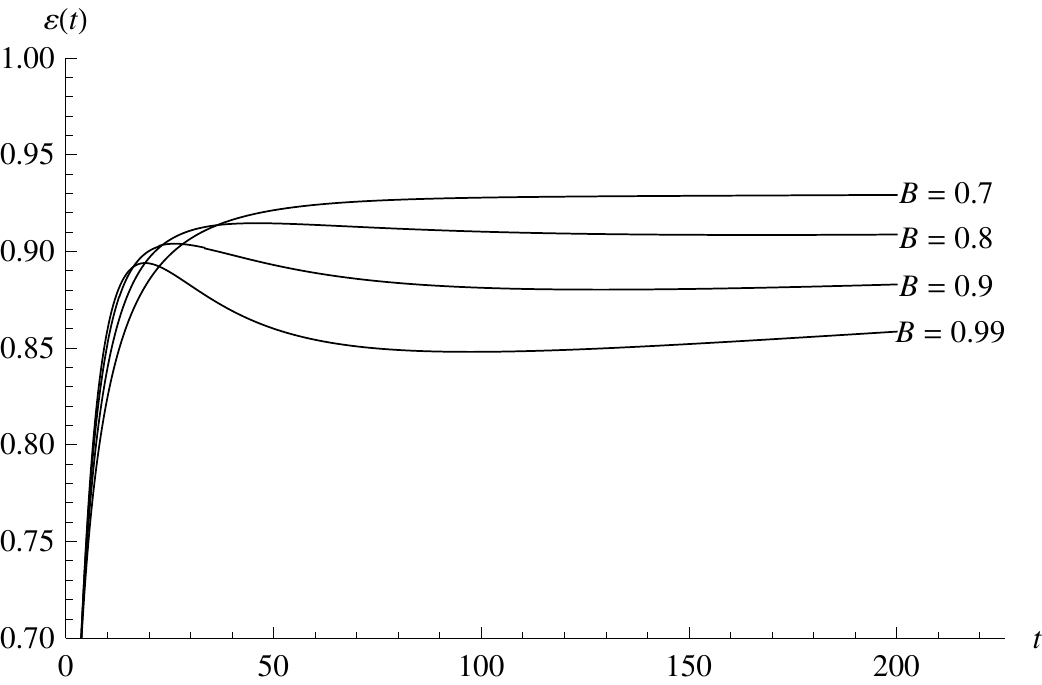}
 \caption{Creep curves for $\protect\alpha =0.4$ and $B\in \{0.7,0.8,0.9,0.99\}$,
 $a=0.8$, $b=0.1$, $t\in [0,200]$.}
 \label{fig-creep-2}
 \end{figure}

 Suppose that $\sigma(t) =H(t)$, i.e.,
 \beq \label{sigma=H}
 H(t) = (1+a\,{}_{0}D_{t}^{\alpha} + 2b\,{}_{0}\bar{D}_{t}^{\beta}) \eps (t),
 \quad t\geq 0.
 \eeq
 By \eqref{n2}, we obtain
 $$
 \eps (t) \approx \frac{H(t) +\sum_{p=1}^{N} \Big\{
 \frac{C_{p-1}(\alpha)}{t^{p+\alpha}}
 +
 \frac{C_{p-1}(\alpha +iB)}{t^{p+\alpha + iB}}
 +\frac{C_{p-1}(\alpha -iB)}{t^{p+\alpha -iB}}
 \Big\}
 V_{p-1}(\eps)(t)}{1+
 \frac{a\,A(N,\alpha)}{t^{\alpha}}+2b\,\Big(
 \frac{A(N,\alpha +iB)}{t^{\alpha +iB}}+\frac{A(N,\alpha -iB)}{t^{\alpha -iB}}
 \Big)},
 $$
 or
 \begin{equation} \label{nn6}
 \eps (t) \approx \frac{H(t) t^{\alpha} + \sum_{p=1}^{N}\Big\{
 \frac{C_{p-1} (\alpha)}{t^{p}}
 +
 \frac{C_{p-1}(\alpha +iB)}{t^{p+iB}}
 +\frac{C_{p-1}(\alpha -iB)}{t^{p-iB}}
 \Big\}
 V_{p-1}(\eps)(t)}{t^{\alpha}+
 a\,A(N,\alpha )+2b\,\Big(
 \frac{A(N,\alpha +iB)}{t^{iB}}+\frac{A(N,\alpha -iB)}{t^{-iB}}
 \Big)}.
 \end{equation}
 By using \eqref{nn6} in \eqref{n3} we obtain
 \beas
 V_{p-1}^{(1)}(\eps)(t) &\approx& t^{p-1} \frac{H(t) t^{\alpha}
 +\sum_{p=1}^{N} \Big\{
 \frac{C_{p-1}(\alpha)}{t^{p}}+
 \frac{C_{p-1}(\alpha +iB)}{t^{p+ iB}}+\frac{C_{p-1}(\alpha -iB)}{t^{p-iB}}
 \Big\}
 V_{p-1}(\eps)(t)}{t^{\alpha} +
 a\,A(N,\alpha)+2b\, \Big(
 \frac{A(N,\alpha +iB)}{t^{iB}} + \frac{A(N,\alpha -iB)}{t^{-iB }}
 \Big)},  \nonumber \\
 V_{p-1}(\eps)(0) &=& 0,
 \quad p=1,2,3,\ldots
 \eeas
 Equation \eqref{sigma=H} may also be solved by contour integration
 \begin{equation} \label{eps-1}
 \eps(t) = K(t)\ast H(t), \quad t\geq 0,
 \end{equation}
 where $K$ is given by \eqref{k}, see \eqref{eps}. Finally, the values of
 $\eps$, at discrete points, could be determined directly from
 $\tilde{\eps}(s) = \frac{1}{s} \tilde{K}(s)$, $\rep s>0$, see \eqref{eps-k-tilda}, by the
 use of Post inversion formula, see \cite{Cohen}.
 Thus,
 \begin{equation} \label{n8}
 \eps(t) = \lim_{n\to \infty} \frac{(-1)^{n}}{n!}
 \bigg[s^{n+1} \frac{d^{n}}{ds^{n}} \frac{1}{s(1+as^{\alpha}+b(s^{\beta}+s^{\bar{\beta}}))}\bigg]_{s=\frac{n}{t}},
 \quad t\geq 0.
 \end{equation}
 In Fig. \ref{fig-creep-3}, \ref{fig-creep-1} and \ref{fig-creep-2}
 we show $\eps$ for several values of
 parameters determined from \eqref{eps-1}. In Fig. \ref{fig-creep-3}, for specified
 values of $t$ we present values of $\eps$, determined from \eqref{nn6},
 with $N=7$, denoted by dots, as well as the values of $\eps$,
 determined by \eqref{n8}, with $n=25$, denoted by squares. As could be seen
 the agreement of results determined by different methods is significant.

 From Fig. \ref{fig-creep-1} and \ref{fig-creep-2} one sees that,
 regardless of the value of $B$, creep curves tend to $\eps=1$. In
 Fig. \ref{fig-creep-1} the creep curves are monotonically increasing,
 while in Fig. \ref{fig-creep-2} creep curves has oscillatory character,
 characteristic for the case when the mass of the rod is not neglected.
 Note that in all cases of $B$, the restriction determined by the dissipation
 inequality is satisfied.

%%%%%%%%%%%%%%%%%%%%%%%%%%%%%%%%%%%%%%%%%%%%%%%%%%
\section{Conclusion}
%%%%%%%%%%%%%%%%%%%%%%%%%%%%%%%%%%%%%%%%%%%%%%%%%%

 In this work, we proposed a new constitutive equation with fractional
 derivatives of complex order for viscoelastic body of the Kelvin-Voigt type.
 The use of fractional derivatives of complex order, together with
 restrictions following from the Second Law of Thermodynamics, represent the
 main novelty of our work. Our results can be summarized as follows.

 \begin{enumerate}
 \item In order to obtain real stress for given real strain, we used two
 fractional derivatives of complex order that
 are complex conjugated numbers, see Theorem \ref{th:sufc}.
 \item The restrictions that follow from the Second Law of Thermodynamics
 for isothermal process implied that the constitutive
 equation must additionally contain a fractional derivative of real order.
 Thus, the simplest constitutive equation that gives real stress for real strain and
 satisfies the dissipativity condition is given by \eqref{nkj}.
 \item We provided a complete analysis of solvability of
 the complex order fractional Kelvin-Voigt model given by
 \eqref{nkj} (see Theorems \ref{th:zeros} and \ref{th:lst}).
 \item We studied stress relaxation and creep problems through equation \eqref{nkj}.
 An increase of $B$ implied that the stress relaxation decreases more rapidly to
 the limiting value of stress, i.e., $\lim_{t\to \infty} \sigma (t) =1$.
 \item We presented numerical experiments when the dissipation inequality is
 satisfied. The creep experiment showed that the increase in the imaginary part of the
 complex derivative $B$ changes the character of creep curve from monotonic
 to oscillatory form, see Fig. \ref{fig-creep-1} and \ref{fig-creep-2}.
 However, the creep curves never cross the value equal to 1. The creep curve
 resembles the form of creep curve when either
 the mass of the rod, see \cite[p.\ 124]{APSZ-book-vol1},
 or inertia of the rheometer, see \cite{JaishankarMcKinley-2012},
 is taken into account.
 Recently, \cite{Zingales} presented experimental results of creep curves
 for some biological materials that exhibit nonmonotonic creep
 curve. Thus, our model \eqref{nkj} is applicable for such a case.
 \item The parameters in the proposed model \eqref{nkj} could be
 determined from experimental data using for instance the least
 square method.
 \item Our further study will be directed to the problems of vibration and wave
 propagation of a rod with finite mass and constitutive equation \eqref{nkj},
 along the lines presented in \cite{APSZ-book-vol1}, see also \cite{Hanyga-qua}.
 \end{enumerate}

%%%%%%%%%%%%%%%%%%%%%%%%%%%%%%%%%%%%%%%%%%%%%%%%%%
 \section*{Acknowledgements}
%%%%%%%%%%%%%%%%%%%%%%%%%%%%%%%%%%%%%%%%%%%%%%%%%%

 We would like to thank Marko Janev for several helpful discussions on the subject.

 This work is supported by
 Projects 174005 and 174024 of the Serbian Ministry of Science,
 and 114-451-1084 of the Provincial Secretariat for Science.

%%%%%%%%%%%%%%%%%%%%%%%%%%%%%%%%%%%%%%%%%%%%%%%%%%
%%%%%%%%%%%%%%%%%%%%%%%%%%%%%%%%%%%%%%%%%%%%%%%%%%

%%%%%%%%%%%%%%%%%%%%%%%%%%%%%%%%%%%%%%%%%%%%%%%%%%
%%%%%%%%%%%%%%%%%%%%%%%%%%%%%%%%%%%%%%%%%%%%%%%%%%
 \end{document}